\documentclass[letterpaper,12pt]{article}
\usepackage{basic_template}
\usepackage{bbm}
\usepackage[numbers]{natbib}
\usepackage{booktabs}
\newcommand{\abs}[1]{\ensuremath{\left| #1 \right|}}

\numberwithin{equation}{section}

\AtBeginEnvironment{lemma}{\crefalias{theorem}{lemma}}
\AtBeginEnvironment{corollary}{\crefalias{theorem}{corollary}}
\AtBeginEnvironment{proposition}{\crefalias{theorem}{proposition}}
\AtBeginEnvironment{assumption}{\crefalias{theorem}{assumption}}
\AtBeginEnvironment{definition}{\crefalias{theorem}{definition}}
\AtBeginEnvironment{example}{\crefalias{theorem}{example}}
\AtBeginEnvironment{remark}{\crefalias{theorem}{remark}}
\AtBeginEnvironment{conjecture}{\crefalias{theorem}{conjecture}}
\AtBeginEnvironment{problem}{\crefalias{theorem}{problem}}

\newcommand\numberthis{\addtocounter{equation}{1}\tag{\theequation}}
\definecolor{NavyBlue}{cmyk}{0.94,0.54,0,0}
\newcommand{\blue}[1]{{\color{NavyBlue}#1}}

\newcommand{\argmax}{\operatorname*{arg\, max}}
\DeclareMathOperator*{\E}{\mathbb{E}}
\newcommand{\R}{\mathbb{R}}
\newcommand{\Proj}{\mathbf{Proj}}

\newcommand{\cS}{\mathcal{S}}
\newcommand{\cA}{\mathcal{A}}

\newcommand{\cT}{\mathcal{T}}
\newcommand{\cP}{\mathcal{P}}
\newcommand{\bH}{\mathbf{H}}

\newcommand{\dS}{{\vert\cS\vert}}
\newcommand{\dA}{{\vert\cA\vert}}

\newcommand{\piref}{{\pi_{\mathrm{ref}}}}

\newcommand{\pitheta}{{\pi_{\theta}}}
\DeclareMathOperator*{\diag}{diag}
\newcommand{\norm}[1]{\left\|#1\right\|}

\begin{document}

\title{Non-Asymptotic Global Convergence of PPO-Clip}
\author{
    Qiming Dai\thanks{School of Mathematical Sciences, Peking University, Beijing 100871, China (e-mail: qmdai25@stu.pku.edu.cn).} \and
    Yin Liu\thanks{Beijing International Center for Mathematical Research, Peking University, Beijing 100871, China (e-mail: yinliu@pku.edu.cn).}\textsuperscript{\phantom{2},*} \and
    Junyu Zhang\thanks{Department of Industrial Systems Engineering and Management, National University of Singapore, Singapore 117576, Singapore (e-mail: junyuz@nus.edu.sg).} \and
    Zaiwen Wen\thanks{Beijing International Center for Mathematical Research and Center for Machine Learning Research, Peking University, Beijing 100871, China (e-mail: wenzw@pku.edu.cn).}
}
\date{September 12, 2026}
\begingroup
\makeatletter
\let\@fnsymbol\@arabic
\makeatother
\maketitle
\endgroup

\begin{abstract}
    Reinforcement learning  has gained attention for modern Large Language Model post-training. The actor-only variants of Proximal Policy Optimization (PPO) are widely applied for their efficiency. These algorithms incorporate a clipping mechanism to improve stability. Besides, a regularization term, such as the reverse KL-divergence or a more general \(f\)-divergence, is introduced to control excessive deviation from a reference policy. Despite their empirical success, a rigorous theoretical understanding of the problem and the algorithm's properties is limited. This paper advances the theoretical foundations of the PPO-Clip algorithm by analyzing a deterministic actor-only PPO algorithm within the general RL setting with \(f\)-divergence regularization under the softmax policy parameterization. We derive a non-uniform Lipschitz smoothness condition and a Łojasiewicz inequality for the considered problem. Based on these properties, we establish non-asymptotic global linear convergence in value gap for the forward KL regularizer. For the reverse KL regularizer, we derive global linear convergence from any finite softmax initialization in both value gap and squared policy distance.
\end{abstract}


\section{Introduction}
Reinforcement Learning (RL) has become a key component of modern Large Language Model (LLM) post-training, where pretrained models are optimized according to reward signals that encode task performance or human preferences. In practice, many RL-based post-training algorithms rely on PPO-style objectives that constrain policy updates through probability-ratio clipping. Although such clipping mechanisms have proved effective in stabilizing large-scale training, the convergence behavior of the resulting clipped policy optimization remains underexplored. This motivates us to study the optimization dynamics induced by clipped policy objectives.

The commonly applied RL algorithms in LLMs are built on the foundations of policy gradient (PG) methods \cite{williams1992simple}. As PG methods are often unstable in practice due to the inexactness of gradient estimation, the Proximal Policy Optimization (PPO) algorithm \cite{schulman2017proximal} is introduced to improve the training performance. It is a practical alternative to Trust Region Policy Optimization (TRPO) \cite{schulman2015trust} that restricts policy updates to a trust region to ensure stability. As an actor-critic method, original PPO requires maintaining both a policy and a value model, which is expensive in both computation and storage, particularly for modern large-scale LLMs. This has motivated research into more efficient actor-only algorithms, such as ReMax \cite{Li2024ReMax}, Group Relative Policy Optimization (GRPO) \cite{shao2024deepseekmath}, and Decoupled Clip and Dynamic sAmpling Policy Optimization (DAPO) \cite{yu2025dapo}, which eliminate the expensive critic model. Nevertheless, these methods still retain PPO's core stability technique: clipping the policy ratio to penalize overly large updates. These methods share the clipped probability-ratio structure, leading us to analyze the mechanism behind it.

Another challenge in LLM post-training is excessive deviation from the reference policy, which may lead to reward overoptimization and degraded generation quality. The training process can also suffer from entropy collapse, where the policy quickly becomes nearly deterministic, limiting exploration and further optimization progress \cite{yu2025dapo}. A common approach to controlling policy deviation is to regularize the objective with a divergence term that penalizes departures from a fixed reference policy. Among such choices, the Kullback--Leibler (KL) divergence has been widely adopted.

However, the standard reverse KL-divergence regularizer can be restrictive. It has been associated with mode-seeking behavior and reduced output diversity. This may cause the policy to concentrate on a narrow set of replies \cite{khalifa2021a, glaese2022improving, wiher2022decoding, perez2022red}. To address this limitation, recent literature has considered the more general class of \(f\)-divergence regularizers \cite{Go2023Aligning, wang2024beyond}, which provides broader regularization geometries. For example, the \(\chi^2\)-divergence, which also belongs to the \(f\)-divergence family, has been applied to improve robustness to reward noise \cite{huang2025correcting}.

\textbf{Contributions.}
Motivated by the above considerations, we study the global convergence behavior of actor-only PPO-Clip with \(f\)-divergence regularization in a finite discounted Markov Decision Process (MDP) under tabular softmax policies. Despite its controlled setting, this problem retains several central difficulties of policy optimization. It therefore provides a tractable framework for characterizing the global convergence behavior of clipped policy optimization. Our main contributions are as follows.

\noindent (1) The \(f\)-divergence regularized value function is formally defined, and its fundamental structural properties under the softmax policy parameterization are analyzed.

\noindent (2) Key properties of the defined objective function are established, including the non-uniform Lipschitz smoothness condition and the Łojasiewicz inequality, under moderate assumptions on the \(f\)-function.

\noindent (3) Building on these properties, the deterministic actor-only PPO-Clip algorithm with forward and reverse KL-divergence regularizers is analyzed. For forward KL, a non-asymptotic global linear convergence rate is established under suitable initialization and stepsize. For reverse KL, global linear convergence in both value gap and squared policy distance is established from any finite softmax initialization.

\textbf{Related Work.}  
PG is one of the most popular algorithms for RL \cite{williams1992simple, sutton1999policy, konda1999actor, kakade2001natural}.
Despite the non-convex landscape of the RL optimization problem \cite{agarwal2021pg}, the convergence of exact policy gradient methods has been extensively studied.
In general policy parameterizations, convergence to stationary points can be established.
Under more structured settings, global convergence to optimal policy has been proven.
For instance, global convergence has been established for the projected PG under direct parameterization \cite{agarwal2021pg, zhang2020variational, xiao2022convergence}, policy mirror descent \cite{xiao2022convergence, lan2023policy, zhan2023policy}, and natural policy gradient (NPG) methods \cite{kakade2001natural, agarwal2021pg}.
Notably, broader structural conditions for PG have been identified that ensure gradient domination, thereby guaranteeing convergence to a global optimum beyond the tabular policy \cite{bhandari2024global}.
Beyond exact gradient methods, a growing number of works have analyzed stochastic PG/NPG algorithms.
Under general policy parameterizations, convergence to stationary points has been established in \cite{xiong2021non, gargiani2022page, paczolay2024sample, mondal2024improved}.
When the policy parameterization satisfies the Fisher non-degeneracy condition or a local bijection condition, sample complexity bounds for converging to a global optimum have been derived in \cite{liu2020improved, fatkhullin2023stochastic, ding2022global,zhang2021convergence}.

When considering regularized RL problems, the global convergence of PG with log-barrier regularization has been analyzed in \cite{williams1992simple, agarwal2021pg, zhang2021sample}.
The convergence properties of PG under the ABC assumption (a relaxed bounded variance assumption) are further investigated in \cite{yuan2022general}.
However, the convergence analysis in these works, despite using a log-barrier regularizer, aims at solving the unregularized objective, leaving the convergence properties of the regularized objective unexplored.
Entropy regularization is another setting widely adopted in practice that encourages exploration and avoids early convergence to local sub-optima \cite{williams1991function, peters2010relative, mnih2016asynchronous, duan2016benchmarking, haarnoja2017reinforcement, nachum2017bridging, neu2017unified, geist2019theory, hazan2019provably, xiao2019maximum, vieillard2020leverage}.
The convergence properties of the entropy-regularized problem for both deterministic and stochastic algorithms have also been investigated in many works \cite{mei2020global, cen2022fast, ding2025twophase}.

Convergence analyses for the PPO algorithm remain relatively limited. Among the two primary PPO objective structures, PPO-KL possesses a simpler analytical form, and its convergence results have been studied in \cite{liu2019neural,holzleitner2021convergence}. For the more widely implemented PPO-Clip objective, theoretical results are scarce due to nondifferentiability. Specifically, Jin et al. \cite{jin2023stationary} have established the stationary-point convergence of PPO-Clip under general policy parameterization, while Huang et al. \cite{huang2024ppo} have analyzed the global convergence of PPO-Clip within a mirror descent framework. However, these results have limitations: the analysis in \cite{huang2024ppo} involves a substantial number of technical constants, making the results less interpretable; while \cite{jin2023stationary} only studies a simplified variant of PPO-Clip. Moreover,  existing convergence analyses largely focus on unregularized problems and do not characterize how divergence regularization affects its global optimization behavior.

\textbf{Notation.} Throughout this paper, \(\|\cdot\|_2\) denotes the  \(l_2\)-norm for vectors and the spectral norm for matrices. We denote the span norm for a vector $x$ as \(\|x\|_{\mathrm{span}} := \max_{i} x_i - \min_{i} x_i\).
\(\delta_a\in \mathbb{R}^{\abs{\cA}}\) denotes the vector in which only the element corresponding to action \(a\) is 1, and \(\delta_s\in \mathbb{R}^{\dS \dA}\) denotes the vector in which only the elements corresponding to state \(s\) are 1. \(\mathbf{1}\) is the all-ones vector, and \(\mathbbm{1}\) is the indicator function. \(\mathbf{I}\) is the identity matrix. For a matrix $A$, the matrix infinity norm is defined as $\|A\|_\infty:=\max_{i}\sum_{j}\abs{A_{ij}}$, which is the maximum absolute row sum.
We denote  \(\theta_s :=\theta(s,\cdot)\), \(\pitheta(s):=\pitheta(\cdot|s)\),     and \(\piref(s):= \piref(\cdot|s)\).
We also denote \(\bH(\theta_s):= \diag(\pitheta(s)) - \pitheta(s) \pitheta(s)^\top\).
For any algorithmic vector \(x_{n,k}\in\mathbb R^{\dS\dA}\), we write \(x_{n,k}(s,\cdot)\) for its state block and \(x_{n,k}(s,a)\) for its action coordinate.

\section{Preliminary}
In this paper, we study the reinforcement learning problem within the framework of an infinite-horizon discounted Markov Decision Process (MDP), defined by the tuple \((\cS,\cA,\cP,r,\gamma)\). Here \(\cS\) is the set of states, and \(\cA\) is the set of actions. The function \(\cP: \cS \times \cA \times \cS \rightarrow[0,1]\) is the state transition probability function, with \(\cP(s'| s,a)\) denoting the probability of transitioning from state \(s\) to state \(s'\) upon taking action \(a\). The reward function is \(r:\cS \times \cA \rightarrow \mathbb{R}\), and \(\gamma \in [0,1)\) is the discount factor. The policy \(\pi:\cS \rightarrow\Delta_{\dA}\) determines the agent's behavior, where \(\pi(a| s)\) is the probability of selecting action \(a\) in state \(s\). In this study, we focus on the softmax policy, defined as \(\pitheta(a| s)=\frac{\exp(\theta_{sa})}{\sum_{a'\in\mathcal A}\exp(\theta_{sa'})}\),
with \(\theta \in \mathbb{R}^{\dS \dA}\) being the policy parameter.

Given an initial state \(s_0\), the interaction between the agent and the MDP generates a trajectory \(\tau=\{s_0,a_0,s_1,a_1,\cdots\}\). 
We use \(\tau\sim \pitheta\) to denote \(\tau\) that follows the probability distribution \(\Pr(\tau | s_0,\pi_\theta)\).
In practice, the initial state also follows a probability distribution \(\rho\), which is typically unknown. Instead, we apply a different distribution \(u\) when calculating the policy gradient. Another important probability distribution is the discounted state distribution \(d^\pitheta_u(s)\) (and \(d^\pitheta_\rho(s)\)), which is the time-averaged discounted probability of visiting state \(s\) given the initial state distribution \(u\) (and \(\rho\)),
\begin{align*}
    d^\pitheta_u(s) = & \E_{s_0\sim u}\Big[(1-\gamma) \sum_{t=0}^\infty \gamma^t\Pr(s_t = s | s_0,\pitheta, \mathcal{P}) \Big].
\end{align*}

\begin{definition}[\(f\)-divergence]\label{def:f-div}
    Let the generator function \(f:(0, \infty) \rightarrow \mathbb{R}\) be a convex function with \(f(1)=0\), the \(f\)-divergence for two distributions \(p\) and \(q\) is defined as
    \[
        D_f(p,q) :=\E_{q(x)}[ f\big(p(x)/q(x)\big) ].
    \]
\end{definition}
The \(f\)-divergence is a general measure of the difference between two distributions. As listed in \cref{tab:f-c-value}, by selecting a different generator function, one can derive various commonly used measures, including forward and reverse KL-divergence, Jensen-Shannon (JS) divergence, \(\chi^2\)-divergence, etc.

\begin{algorithm}[hbt]
    \caption{PPO-Clip}\label{alg:PPO-Clip}
    \begin{algorithmic}[1]
        \REQUIRE Initialization: \(\pi_{\theta_{1,1}} = \piref\), stepsize \(\{\eta_{n,k}\}\)
        \FOR{\(n=1,\cdots, N\)}
        \FOR{\(k=1,\cdots,K\)}
        \STATE \(\theta_{n,k+1} =  \theta_{n,k} + \eta_{n,k} \nabla_{\theta_{n,k}} \mathcal{L}_n(\theta_{n,k})\)
        \ENDFOR
        \STATE \(\theta_{n+1,1}=\theta_{n,K+1}\)
        \ENDFOR
    \end{algorithmic}
\end{algorithm}

PPO-style algorithms used in LLM post-training, including PPO, GRPO, and DAPO, optimize surrogate objectives that approximate the original policy objective. These methods share several high-level ingredients with \cref{alg:PPO-Clip}, including a fixed old policy, a clipped probability-ratio surrogate, and multiple inner updates, while differing in sampling, advantage estimation, regularization, and policy parameterization. \Cref{alg:PPO-Clip} has a double-loop structure. At the \(n\)-th outer iteration, the old policy \(\pi_{n,1}\) is fixed,  which determines the sampling distribution. The algorithm then performs several gradient ascent steps on the corresponding surrogate objective in the inner loop.
The operator \({\rm clip}(v,1-\varepsilon_l,1+\varepsilon_h)\) clips \(v\) to the interval \([1-\varepsilon_l,1+\varepsilon_h]\). Throughout this paper, we allow asymmetric clipping thresholds \(\varepsilon_l\) and \(\varepsilon_h\) for generality, as in DAPO, whereas GRPO typically adopts the symmetric choice \(\varepsilon_l=\varepsilon_h\).

We adopt the shorthand \(\E_{\pi_{n,1}} := \frac{1}{1-\gamma}\E_{s\sim d^{\pi_{n,1}}_u, a\sim \pi_{n,1}(s)}\) throughout the rest of the paper. Using the notation 
\[r_{\theta}^n(s,a):=\frac{\pitheta(a| s)}{\pi_{n,1}(a| s)}, \qquad  r_{\rm ref}^n(s,a):= \frac{\piref(a |s)}{\pi_{n,1}(a|s)},\]
the surrogate problem with a general \(f\)-divergence regularizer can be written in the following unified form:
\begin{align*}
     \mathcal{L}_{n}(\theta) &:= \E_{\pi_{n,1}}  \Big[\min\big\{r_\theta^n(s,a) \tilde{A}^{\pi_{n,1}}_\lambda(s,a), \\
&\qquad \text{clip}\big(r_\theta^n(s,a);1-\varepsilon_l, 1+\varepsilon_h\big) \tilde{A}^{\pi_{n,1}}_\lambda(s,a)\big\} \\
&\qquad -  \lambda D_f(\pitheta(s),\piref(s)) \Big],
\end{align*}
where \(\tilde{A}_\lambda^{\pi_{n,1}}\) denotes the regularized advantage function whose rigorous definition is given in \cref{lem:gradient}.

\section{Properties of \texorpdfstring{\(f\)}{f}-divergence regularized value function}

As motivated in the introduction, the problem we consider is the \(f\)-divergence regularized value function with the form
\begin{align*}
    \tilde{V}^{\pi_\theta}_\lambda(\rho) = \! \E_{\rho,\pi_\theta}\! \Big[ &\sum_{t=0}^{\infty} \!\gamma^t r(s_t,a_t) \! -\! \lambda \sum_{t=0}^{\infty} \!\gamma^t D_f\big(\pi_\theta (s_t), \pi_{\text{ref}}(s_t)\big) \Big].
\end{align*}
Such an objective function is a generalization of the widely studied entropy-regularized reward function \cite{williams1991function,mnih2016asynchronous,mei2020global}. As \(\piref\) is a fixed reference policy, we can simplify the notation of \(f\)-divergence as \(D_f(\theta_s)\). This function can be separated into  \(\tilde{V}^\pitheta_\lambda(\rho) = V^\pitheta(\rho)  - \lambda V^\pitheta_f(\rho)\), where the two value functions are associated with  reward and \(f\)-divergence, respectively:
\begin{align*}
    V^\pitheta(\rho)    := \E_{ \rho, \pitheta}\Big[\sum_{t=0}^\infty \gamma^t r(s_t,a_t) \Big ], \qquad V^\pitheta_f(\rho)  := \E_{ \rho,  \pitheta} \Big[ \sum_{t=0}^\infty\gamma^t D_f(\theta_{s_t}) \Big].
\end{align*}
Similarly, the regularized reward can be defined as
\begin{align*}
    \tilde{r}^\pitheta_\lambda(s_t,a_t) := r(s_t,a_t) - \lambda  D_f(\theta_{s_t}).
\end{align*}
We also make some assumptions about the reference policy and initial state distribution, which are commonly adopted in the theoretical analysis of RL.
\begin{assumption}\label[assumption]{asup:min-piref}
    The reference policy \(\piref\) has a uniform lower bound on each entry, i.e., \(c_{\piref}:=\min_{a,s}\piref(a | s)>0\). The reward $r(s,a)$ is bounded between \([0,r_{\mathrm{max}}]\).
\end{assumption}
Let \(\theta_{\mathrm{ref}}\) be the corresponding parameter of the reference policy \(\piref\). When \(\pi_{\mathrm{ref}}\) is the uniform policy that encourages exploration, \(c_{\piref}=\frac{1}{\dA}\). The reverse KL-regularizer reduces to entropy, the forward KL-regularizer reduces to log-barrier.
\begin{assumption}\label[assumption]{asup:min-u}(Explorative start)
    The starting state distribution \(u\) is strictly positive, i.e., \(u(s)\geq c_u >0,\forall s\in\cS\).
\end{assumption}
This assumption implies \(d_u^\pi(s)\ge (1-\gamma)u(s)>0\) for all \(s\in\cS\) and all \(\pi\in\Pi\). Depending on the separation of the value function, the gradient can be calculated separately in the form of \(Q\)-function and the advantage function individually.
\begin{lemma}\label{lem:gradient}
    Let \(w^{\pitheta}_s \in \mathbb{R}^{\abs{\cA}}\) be a vector with \([w^{\pitheta}_s]_a = w^{\pitheta}_{sa}:= \frac{\pitheta(a | s)}{\piref(a | s)}\). The gradient of \(\tilde V_\lambda^\pitheta (u)\) with respect to \(\theta\) is
    \begin{align*}
        &\frac{\partial \tilde{V}^\pitheta_\lambda(u)}{\partial \theta_{sa}} =  \frac{d^\pitheta_u(s)}{1-\gamma} \pitheta(a | s) \Big[\tilde{A}_\lambda^\pitheta(s,a) -\lambda \Big( f' (w^{\pitheta}_{sa} )         - \sum_{a'}\pitheta(a' | s) f' \left(w^{\pitheta}_{sa'}\right) \Big)\Big],
    \end{align*}
    where the total advantage function of the reward and regularizer is defined as
    \begin{align*}
        \tilde{A}^\pitheta_\lambda(s,a) =& \left( Q^\pitheta(s,a)- \lambda Q^\pitheta_f(s,a) \right)-  \left(V^\pitheta(s) - \lambda V^\pitheta_f(s)\right).
    \end{align*}
    \end{lemma}
    \begin{proof}
    The Q-function of \(V^\pitheta(u)\) is 
    \begin{align*}
       Q^\pitheta (s,a)= r(s,a) + \gamma \sum_{s'} \cP(s' | s,a) V^\pitheta(s').
    \end{align*}
    follows the standard policy-gradient derivation, so we highlight only the key steps. The gradient of \(V^\pitheta_f(u)\) is obtained analogously. By the definition, we know
    \begin{align*}
        V^\pitheta_f(u) = \E_{s\sim u}\Big[\sum_a \pitheta(a | s) Q^\pitheta_f(s,a) \Big], \quad \text{with } Q^\pitheta_f(s,a) = D_f(\theta_s)+ \gamma \sum_{s'}\cP(s' | s,a) V^\pitheta_f(s').
    \end{align*}
    Taking gradient with respect to \(\theta\), we get
    \begin{align*}
        \nabla_\theta  V^\pitheta_f(u)  & =  \underset{s \sim u}{\E}\Big[\sum_a  \nabla_\theta \pitheta(a | s)Q^\pitheta_f(s,a) + \sum_a \pitheta(a | s) \nabla_\theta Q^\pitheta_f(s,a)   \Big]         \\
        &=  \underset{s \sim u}{\E} \Big[\sum_a \nabla_\theta \pitheta(a | s)\left(Q^\pitheta_f(s,a)+ f'\left(w^{\pitheta}_{sa} \right) \right) \\
        &\qquad + \gamma \sum_a \pitheta(a | s) \sum_{s'}\cP(s' | s,a) \nabla_\theta V^\pitheta_f(s')\Big]                                                 \\
        & =  \frac{1}{1-\gamma}\sum_s d^\pitheta_u(s) \sum_a \nabla_\theta \pitheta(a | s)  \left(Q^\pitheta_f(s,a)+ f'\left(w^{\pitheta}_{sa} \right) \right).
    \end{align*}
    The last equality follows the proof of the policy gradient theorem. It is derived by recursively expanding the gradient of the value function and using the definition of the discounted state visitation distribution, \(d^{\pi_\theta}_u(s)\). For a detailed derivation, see Chapter 13 of \cite{sutton1998reinforcement}.
   Combined with the gradient of the value function, we have the total gradient as
    \begin{align*}
        \nabla_\theta \tilde{V}^\pitheta_\lambda(u)  = &  \frac{1}{1-\gamma}\sum_s d^\pitheta_u(s) \sum_a \nabla_\theta \pitheta(a | s)  \left(Q^\pitheta(s,a) - \lambda Q^\pitheta_f(s,a) - \lambda f'\left( w^{\pitheta}_{sa}\right) \right).
    \end{align*}
    For \(s'\neq s\), it holds \(\frac{\partial \pitheta(a | s)}{\partial \theta_{s'}}=0\), hence
    \begin{align*}
        \frac{\partial  V_f^\pitheta(u)}{\partial \theta_s} = & \frac{d^\pitheta_u(s)}{1-\gamma} \sum_{a'} \frac{\partial\pitheta(a' | s)}{\partial\theta_s} \left(Q_f^\pitheta(s,a') + f' \left(w^{\pitheta}_{sa'} \right) \right) \\
        =                                                     & \frac{d^\pitheta_u(s)}{1-\gamma} \mathbf{H}(\theta_s) \left(Q_f^\pitheta(s) +f' \left(w^{\pitheta}_s \right) \right),
    \end{align*}
    where $Q_f^\pitheta(s):=Q_f^\pitheta(s,\cdot)$.
    The whole derivative is
    \begin{align*}
        \frac{\partial \tilde{V}^\pitheta_\lambda(u)}{\partial\theta_s} = \! \frac{d^\pitheta_u(s)}{1-\gamma}\bH(\theta_s)\!\left(Q^\pitheta(s) \!-\! \lambda Q^\pitheta_f(s)\! -\! \lambda f' \left(w^{\pitheta}_s \right) \right).
    \end{align*}
    For the derivative with respect to \(\theta_{sa}\), it holds that \begin{align*}
        \frac{\partial \pitheta(a | s)}{\partial\theta_{sa'}} = (\mathbbm{1}_{a=a'} - \pitheta(a' | s)) \pitheta(a | s),
    \end{align*}
    hence the partial derivative follows
    \begin{align*}
          \frac{\partial  V_f^\pitheta(u)}{\partial \theta_{sa}} =&  \frac{d^\pitheta_u(s)}{1-\gamma} \pitheta(a | s)  \Big[Q_f^\pitheta(s,a ) + f' \left(w^{\pitheta}_{sa} \right) - \sum_{a'} \pitheta(a' | s)\left(Q_f^\pitheta(s,a')+  f' \left(w^{\pitheta}_{sa'} \right) \right) \Big]                                                                          \\
        = & \frac{d^\pitheta_u(s) }{1-\gamma } \pitheta(a | s) \Big[A^\pitheta_f(s,a) + f' \left(w^{\pitheta}_{sa} \right) - \sum_{a'}\pitheta(a' | s) f' \left(w^{\pitheta}_{sa'} \right) \Big].
    \end{align*}
    Combined with the gradient \(\nabla_\theta V^{\pitheta}(u)\), the whole derivative can be obtained as the claim.
    \end{proof}

Next, we will show some generalized properties of the regularized value function that will be used when analyzing the convergence rate of the PPO-Clip algorithm.
\begin{lemma}\label{lem:bounded-psi}
    Let \(\psi_\theta^s(s,a):= \nabla_{\theta_s} \log \pitheta(a | s)\) be the score function. For any \(\theta,\theta'\in \mathbb{R}^{\abs{\cS}\abs{\cA}}\) and any \(s\in\cS, a\in\cA\), the softmax policy parameterization satisfies:
    \begin{align}
         & \norm{\psi_\theta^s(s,a)}_2  \leq \sqrt{2}\label{ineq:regularity-lemma-1}                                                                                       \\
         & \norm{\nabla_{\theta_s} \psi_\theta^s(s,a)}_2 \leq 1\label{ineq:regularity-lemma-2}                                                                            \\
         & \norm{\nabla_{\theta_s} \pitheta(a|s) - \nabla_{\theta_s} \pi_{\theta'}(a | s) }_2 \leq 3 \norm{\theta_s - \theta_s'}_2. \label{ineq:regularity-pitheta-smooth}
    \end{align}
\end{lemma}
\begin{proof}
To show the first inequality, the partial derivative of \(\log \pitheta(a| s)\) is calculated as  \(\psi_\theta^s(s,a) = \delta_{a} - \pitheta(s)\), and  the gradient is bounded as
    \begin{align*}
        \norm{\psi_\theta^s(s,a)}_2^2 = \sum_{a'} \left(\mathbbm{1}_{a=a'} - \pitheta(a' | s)\right)^2  =1-2\pitheta(a | s) + \sum_{a'} \pitheta(a' | s)^2 \leq 2.
    \end{align*}
    For the second inequality, a direct computation yields \(\nabla_{\theta_s} \psi_\theta^s(s,a) = -\bH(\theta_s)\).
     Lemma 22 of \cite{mei2020global} shows the eigenvalues of \(\mathbf{H}(\theta_s)\) lie in \([0,1)\), so \(\norm{\bH(\theta_s)}_2 \leq 1\) and the result is derived. And the next inequality holds automatically as the uniformly bounded Jacobian implies \(\psi_\theta^s(s,a)\) is 1-Lipschitz continuous.

    To show the last inequality, define \(\Delta :=\theta'-\theta\), \(\theta(l) := \theta + l\Delta \). By the Fundamental Theorem of Calculus, we have
    \begin{align*}
          \nabla_{\theta_s} \pitheta(a | s)- \nabla_{\theta_s} \pi_{\theta'}(a | s) &=  \int_0^1 \nabla^2_{\theta_s} \pi_{\theta(l)}(a | s) \Delta_s dl                                                                      \\
         & = \int_0^1 \! \pi_{\theta(l)}(a | s) (  \psi_{\theta(l)}^s (s,a) \psi_{\theta(l)}^s (s,a)^\top - \bH(\theta_s(l))) \Delta_s dl.
    \end{align*}
    Take the norm on both sides, with \eqref{ineq:regularity-lemma-1} and \eqref{ineq:regularity-lemma-2}, we have
    \begin{equation*}
             \norm{\nabla_\theta \pitheta(a | s) - \nabla_\theta \pi_{\theta'}(a | s)}_2 \leq \int_0^1 \norm{\Delta_s}_2\Big( \norm{\bH(\theta_s(l)) }_2  + \norm{\psi^s_{\theta(l)}(s,a)}^2_2\Big)dl  \leq 3 \norm{\Delta_s}_2 . \qedhere
    \end{equation*}
\end{proof}

\begin{lemma}[Performance difference lemma]\label{lem:performance-difference}For two arbitrary policies \(\pi_1\) and \(\pi_2\), it holds that
    \begin{align*}
       \tilde{V}^{\pi_1}_\lambda(\rho) \!- \! \tilde{V}^{\pi_2}_\lambda(\rho) =& \frac{1}{1-\gamma}\sum_s d^{\pi_1}_\rho(s)\!  \Big[\sum_a (\pi_1(a | s) - \pi_2(a | s))  \tilde{Q}^{\pi_2}_\lambda(s,a)\\
       & \quad  - \lambda \big(D_f(\pi_1(s) ,\piref(s))\! -\! D_f(\pi_2(s),\piref(s)) \big)  \Big].
    \end{align*}
\end{lemma}
\begin{proof}
For any given initial state \(s_0\), by the definition of the regularized value function, we have
    \begin{align*}
   & \tilde{V}^{\pi_1}_\lambda(s_0) - \tilde{V}^{\pi_2}_\lambda(s_0) \\
        =                                 & \underset{\tau\sim \pi_1}{\E}\Big[\sum_{t=0}^\infty\gamma^t (\tilde{r}^{\pi_1}_\lambda(s_t,a_t) \!+ \!\tilde{V}^{\pi_2}_\lambda(s_t) - \tilde{V}^{\pi_2}_\lambda(s_t))\Big] \!-\! \tilde{V}^{\pi_2}_\lambda(s_0) \\
        \overset{(i)}{=} & \underset{\tau\sim \pi_1}{\E}\Big[\sum_{t=0}^\infty\gamma^t (\tilde{r}^{\pi_1}_\lambda(s_t,a_t) + \gamma\tilde{V}^{\pi_2}_\lambda(s_{t+1}) - \tilde{V}^{\pi_2}_\lambda(s_t) )\Big]                        \\
        =                                 & \underset{\tau\sim \pi_1}{\E}\Big[\sum_{t=0}^\infty\gamma^t \Big(\tilde{r}^{\pi_2}_\lambda(s_t,a_t) \!+\! \gamma\E[\tilde{V}^{\pi_2}_\lambda(s_{t+1}) | s_t,a_t ]\!-\! \tilde{V}^{\pi_2}_\lambda(s_t)    \\
                                          & \quad -\lambda(D_f(\pi_1(s_t),\piref(s_t)) - D_f(\pi_2(s_t), \piref(s_t))) \Big) \Big]                                                                           \\
        =                                 & \underset{\tau\sim \pi_1}{\E}\Big[\sum_{t=0}^\infty\gamma^t\Big(\tilde{Q}^{\pi_2}_\lambda(s_t,a_t)- \tilde{V}^{\pi_2}_\lambda(s_t)-\lambda(D_f(\pi_1(s_t),\piref(s_t)) - D_f(\pi_2(s_t), \piref(s_t))) \Big)\Big].
    \end{align*}
The equality \((i)\) follows from algebraic rearrangement of the series
\[\sum_{t=0}^\infty \gamma^t \tilde{V}^{\pi_2}_\lambda(s_t) - \tilde{V}^{\pi_2}_\lambda(s_0) = \sum_{t=1}^{\infty}\gamma^t \tilde{V}^{\pi_2}_\lambda(s_t) = \sum_{t=0}^\infty \gamma^{t+1} \tilde{V}^{\pi_2}_\lambda(s_{t+1}).\]     
Taking expectation over the initial distribution \(\rho\), and then converting the expectation to the sum over the discounted state-visitation distribution, we can derive the final result.
\end{proof}

\begin{corollary}[Value sub-optimality]\label{cor:value-sub-optimality} 
Let \(\pi^*\in\Delta_{\cA}^{\cS}\) be a regularized Bellman-optimal policy, i.e., a policy that is optimal from every state. Then, for any policy \(\pi\) and any initial distribution \(\rho\), it holds    
    \begin{align*}
        \tilde{V}^{\pi^*}_\lambda(\rho) - \tilde{V}^\pi_\lambda(\rho) = & \frac{1}{1-\gamma}\sum_s d^{\pi}_\rho(s) \Big[\sum_a (\pi^*(a | s) - \pi(a | s)) \tilde{Q}^{\pi^*}_\lambda(s,a)\\
        & \quad   -\lambda\big(D_f(\pi^*(s),\piref(s)) - D_f(\pi(s), \piref(s))\big) \Big].
    \end{align*}
    Furthermore, for all $s\in \cS$, 
    \begin{align*}
         & \sum_a (\pi^*(a | s) - \pi(a | s)) \tilde{Q}^{\pi^*}_\lambda(s,a)     -\lambda\Big(D_f(\pi^*(s),\piref(s)) - D_f(\pi(s), \piref(s))\Big) \geq 0.
    \end{align*}
\end{corollary}

\begin{proof}
  In \cref{lem:performance-difference}, replacing \(\pi_1\) with \(\pi\), \(\pi_2\) with \(\pi^*\). Multiplying both sides by -1 proves the first equation. 
    To show the non-negativity, we have
    \begin{align*}
             & \sum_a (\pi^*(a | s) - \pi(a | s)) \tilde{Q}^{\pi^*}_\lambda(s,a) -\lambda\Big(D_f(\pi^*(s),\piref(s)) - D_f(\pi(s), \piref(s)) \Big) \\
        & =     \tilde{V}^{\pi^*}_\lambda(s) - \sum_a \pi(a | s)\tilde{Q}^{\pi^*}_\lambda(s,a)  -\lambda D_f(\pi^*(s),\piref(s)) + \lambda D_f(\pi(s), \piref(s)) \\
        & =    \tilde{V}^{\pi^*}_\lambda(s) - \sum_a \pi(a | s) \Big(r(s,a) - \lambda D_f(\pi(s), \piref(s)) + \gamma \sum_{s'}\cP(s' | s,a)\tilde{V}^{\pi^*}_\lambda(s')\Big)\geq  0.
    \end{align*}
    The last inequality holds because, by the Bellman optimality equation, we have
    \begin{equation*}
           \tilde{V}^{\pi^*}_\lambda(s) = \max_{\pi\in\Delta(\cA)} \sum_a \pi(a | s) \Big(r(s,a) - \lambda D_f(\pi(s), \piref(s))+ \gamma \sum_{s'}\cP(s' | s,a)\tilde{V}^{\pi^*}_\lambda(s')\Big). \qedhere
    \end{equation*}
\end{proof}

\subsection{Non-uniform Lipschitz smoothness of the objective function}
In this part, we will derive the smoothness condition of \(\tilde{V}^\pitheta_\lambda(\rho)\) and generalize the result of the entropy-regularized value function in \cite{mei2020global} to \(f\)-divergence regularizer. As Agarwal et al. \cite{agarwal2021pg} have proved the unregularized value function \(V^{\pitheta}(\rho)\) is \( \frac{8r_{\rm max}}{(1-\gamma)^3} \)-smooth. We only need to derive the smoothness condition for the regularizer. Here, we first give the upper bound of the second-order directional derivative of the regularizer, then use this intermediate result to show the smoothness of the objective function. 

Let $F_\theta:\R^k\to \R^d$ be a vector-valued function. Define the first differential applied to a direction $v\in\R^k$ as $\mathcal{D}_\theta F_\theta[v]:=\left(\nabla_\theta F_\theta\right)^\top v$. The bilinear second differential applied twice to $v$ is the vector $\mathcal{D}_\theta^2 F_\theta[v,v]\in \R^d$.

\begin{lemma}[Smoothness Framework]\label{lem:smooth-frame}
    For any \(\theta,v\in \mathbb{R}^{\dS\dA}\), and any state \(s\in\cS\), the following inequality holds
    \begin{align*}
        \left\vert  \mathcal{D}_\theta^2 V_f^{\pi_{\theta}}(s) [v,v]\right\vert &\le   \frac{4 \gamma \|v\|_2}{(1 - \gamma)^2} \norm{ \mathcal{D}_\theta D_f(\theta)[v]}_\infty \\
        &\quad + \frac{8 \gamma \|v\|_2^2}{(1 - \gamma)^3}  \left\| D_f(\theta)\right\|_\infty + \frac{1}{1 - \gamma} \norm{\mathcal{D}_\theta^2 D_f(\theta)[v,v]}_\infty.
    \end{align*}
\end{lemma}
\begin{proof}
    This proof generalizes the main idea of Lemma 14 in \cite{mei2020global} to accommodate a general \(f\)-divergence regularizer. For completeness, we include the full proof below.
    
   Noting that \( \frac{\partial \pi_{\theta}(a | s)}{\partial \theta(s', \cdot)} = 0 \) for \( s \neq s' \), we have
    \begin{equation}\label{eq:pi-alpha-deriv}
        \mathcal{D}_\theta \pi_\theta(a|s)[v]=v_s^\top \mathbf{H}(\theta_s) \delta_a
    \end{equation}
    where \(v_s\) denotes  \(v(s,\cdot)\). Next, we bound the sum of the absolute values of the first-order derivative over actions as
    \begin{align*}
         & \sum_a \left\vert  \mathcal{D}_\theta \pi_\theta(a|s)[v] \right\vert = \sum_a \left\vert  v_s^\top \bH(\theta_s) \delta_a \right\vert                              \\
         \leq & \sum_a \pi_{\theta}(a | s) \left(\abs{v_s^\top \delta_a} + \abs{\pitheta(s)^\top v_s}\right)\le 2\sum_{a}  \pi_{\theta}(a| s)\abs{v_{sa}}   \leq 2 \norm{v_s}_\infty \leq 2  \|v\|_2.
    \end{align*}
    For the second-order derivative along the direction \(v_s\), we have
    \begin{align*}
         & \quad \sum_a \left\vert  \mathcal{D}_\theta^2 \pi_\theta(a|s)[v,v] \right\vert
         = \sum_a \left\vert  v_s^\top \nabla_{\theta_s}^2 \pi_\theta(a|s)v_s\right\vert                   \\
         & = \sum_a \pitheta(a | s) \abs{v_s^\top \Big((\delta_a - \pitheta(s)) (\delta_a - \pitheta(s))^\top - \bH(\theta_s) \Big) v_s}               \\
         & \leq \sum_a \pi_{\theta}(a | s) \Big( \left\vert  v_s^\top \delta_a {\delta_a}^\top v_s \right\vert  + \left\vert  v_s^\top \delta_a \pitheta(s)^\top v_s \right\vert  + \left\vert  v_s^\top \pitheta(s) \delta_a^\top v_s \right\vert  \\
         & \qquad+ 2 \left\vert  v_s^\top \pitheta(s) \pitheta(s)^\top v_s \right\vert + \left\vert  v_s^\top \diag(\pitheta(s)) v_s \right\vert  \Big)   \\
         &\leq 6 \|v_s\|_2^2 \leq 6 \|v\|_2^2.
    \end{align*}
   Let \( P_\theta\) be the transition matrix of \( \pi_{\theta} \) with entries defined as
    \(\big[ P_\theta \big]_{(s,s')} := \sum_a \pi_{\theta}(a | s) \mathcal{P}(s' | s,a),\)
    we can compute its first and second order derivatives as
    \begin{align*}
        \big[ \mathcal{D}_\theta P_\theta[v]\big]_{(s,s')} &= \sum_a \left( \mathcal{D}_\theta \pi_\theta(a|s)[v]\right) \mathcal{P}(s' | s,a), \\
        \left[ \mathcal{D}_\theta^2 P_\theta [v,v] \right]_{(s,s')} &= \sum_a \left( \mathcal{D}_\theta^2 \pi_\theta(a|s) [v,v] \right) \cP(s' | s,a).
    \end{align*}
    For any vector \(\nu \in \mathbb{R}^{\dS} \), we have the following two bounds.
    \begin{align*}
     \left\|  \mathcal{D}_\theta P_\theta [v]\,\nu \right\|_{\infty} &=  \max_s \Big\vert  \sum_{s',a}  \left( \mathcal{D}_\theta \pi_\theta(a|s)[v] \right) \mathcal{P}(s' | s,a) \nu(s') \Big\vert   \\
         &\leq  \max_s \sum_{s',a}  \mathcal{P}(s' | s,a) \left\vert  \mathcal{D}_\theta \pi_\theta(a|s)[v] \right\vert\cdot\|\nu\|_{\infty} \leq  2 \|v\|_2 \|\nu\|_{\infty}. \numberthis \label{smo-lem:p-1}
    \end{align*}
    \begin{align*}
        \left\| \mathcal{D}_\theta^2 P_\theta [v,v] \nu \right\|_{\infty} & =  \max_s \Big\vert  \sum_{s',a}  \left(\mathcal{D}_\theta^2 \pi_\theta(a|s) [v,v]\right) \mathcal{P}(s' | s,a) \nu(s') \Big\vert   \\
        & \leq  \max_s \sum_{s',a} \mathcal{P}(s' | s,a) \left\vert  \mathcal{D}_\theta^2 \pi_\theta(a|s) [v,v] \right\vert\cdot\|\nu\|_{\infty} \leq  6 \|v\|_2^2 \|\nu\|_{\infty}. \numberthis \label{smo-lem:p-2}
    \end{align*}
    Consider the accumulated visitation matrix
    \(M_\theta := \left(\mathbf{I} - \gamma P_\theta\right)^{-1} = \sum_{t=0}^{\infty} \gamma^t P_\theta^t.\)
    This Neumann series converges because \(P_\theta\) is a row-stochastic transition matrix (hence has spectral radius at most 1) and the fact \(\gamma\in (0,1)\) ensuring that the spectral radius of \(\gamma P_\theta\)  is strictly less than 1. Consequently, \(\mathbf{I} - \gamma P_\theta\) is invertible, and the series representation is valid. Its \( (s,s') \)-entry has the form 
    \begin{align*}
        \delta_s^\top M_\theta \delta_{s'} &= \sum_{t=0}^{\infty} \gamma^t \delta_s^\top P_\theta^t \delta_{s'} = \sum_{t=0}^{\infty} \gamma^t \Pr \left\{ s_t = s' | s_0 = s, \pi_{\theta} \right\} = \frac{d_s^{\pi_{\theta}}(s')}{1-\gamma} \geq 0,
    \end{align*}
    which is proportional to the state visitation starting from \( s \). Furthermore, we have that the row-sum satisfies \(  \delta_s^\top M_\theta \mathbf{1} = \frac{1}{1 - \gamma} \), for all \( s \in \mathcal{S} \), which implies that for any \( v_s \in \mathbb{R}^{\dS} \),
    \begin{align*}
        \| M_\theta v_s \|_\infty = \max_s \abs{\delta_s^\top M_\theta v_s} \le \max_{s } \delta_s^\top M_\theta \mathbf{1} \| v_s \|_\infty \le \frac{1}{1 - \gamma} \| v_s \|_\infty. \numberthis  \label{smo-lem:m-1}
    \end{align*}
    Hence, the regularizer is of the form
    \begin{equation*}
        V_f^{\pi_{\theta}}(s) = \sum_{s'}\frac{d_s^{\pi_{\theta}}(s')}{1-\gamma} D_f(\theta_{s'}) = \delta_s^\top M_\theta D_f(\theta),
    \end{equation*}
which gives the vector form  \(V_f^{\pi_{\theta}}  = M_\theta D_f(\theta)\), where $D_f(\theta)\in \mathbb{R}^{\abs{\cS}}$ is the stack of $D_f(\theta_s)$.
    Differentiating the identity \( (\mathbf{I} - \gamma P_\theta)\,M_\theta = \mathbf{I} \) at direction $v$, we obtain
    \[\mathcal{D}_\theta M_\theta[v]=\gamma\, M_\theta\left(\mathcal{D}_\theta P_\theta[v]\right) M_\theta.\]
    Differentiating \( V_f^{\pi_{\theta}}(s) \) with respect to \( l \) gives
    \begin{align*}
        \mathcal{D}_\theta V_f^{\pi_{\theta}}[v]= & \gamma M_\theta \left(\mathcal{D}_\theta P_\theta[v]\right) M_\theta D_f(\theta)  + M_\theta \left(\mathcal{D}_\theta D_f(\theta)[v]\right).
    \end{align*}
    For the second-order derivative, we have
    \begin{align*}
           \mathcal{D}_\theta^2 V_f^{\pi_{\theta}}[v,v] &= 2\gamma^2 M_\theta \left(\mathcal{D}_\theta P_\theta [v]\right) M_\theta \left(\mathcal{D}_\theta P_\theta [v]\right) M_\theta D_f(\theta) + \gamma M_\theta \left(\mathcal{D}_\theta^2 P_\theta[v,v]\right) M_\theta D_f(\theta) \\
         & \quad + 2\gamma M_\theta \left(\mathcal{D}_\theta P_\theta [v]\right) M_\theta \left(\mathcal{D}_\theta D_f(\theta) [v]\right) + M_\theta \left(\mathcal{D}_\theta^2 D_f(\theta)[v,v]\right).
    \end{align*}
    Next, we will bound the above four terms individually at \(l=0\). For the first term, it holds
    \begin{align*}
         & \norm{ M_\theta \left(\mathcal{D}_\theta P_\theta [v]\right) M_\theta \left(\mathcal{D}_\theta P_\theta [v]\right) M_\theta D_f(\theta)}_\infty \le \frac{4 \|v\|_2^2}{(1 - \gamma)^3} \cdot \| D_f(\theta) \|_\infty,
    \end{align*}
    where the last inequality follows from \eqref{smo-lem:p-1} and \eqref{smo-lem:m-1}. For the second term, by utilizing \eqref{smo-lem:p-2} and \eqref{smo-lem:m-1}, we have
    \begin{align*}
        & \norm{ M_\theta \left(\mathcal{D}_\theta^2 P_\theta[v,v]\right) M_\theta D_f(\theta)}_\infty  \! \le \! \frac{6 \|v\|_2^2}{(1 - \gamma)^2}  \| D_f(\theta) \|_\infty.
    \end{align*}
    For the third term, we have
    \begin{align*}
        & \norm{ M_\theta \left(\mathcal{D}_\theta P_\theta [v]\right) M_\theta \left(\mathcal{D}_\theta D_f(\theta) [v]\right)}_\infty \le   \frac{2 \|v\|_2}{(1 - \gamma)^2} \left\|\mathcal{D}_\theta D_f(\theta) [v]\right\|_\infty.
    \end{align*}
    Finally, for the last term, we can obtain
    \begin{align*}
        & \norm{M_\theta \left(\mathcal{D}_\theta^2 D_f(\theta)[v,v]\right)}_\infty \le \frac{1}{1 - \gamma}  \left\| \mathcal{D}_\theta^2 D_f(\theta)[v,v] \right\|_\infty.
    \end{align*}
    Combining these bounds and noting \(\frac{8 \gamma^2 }{(1 - \gamma)^3} + \frac{6\gamma}{(1-\gamma)^2}\leq \frac{8\gamma}{(1-\gamma)^3}\) for \(\gamma\in(0,1)\) yields the desired inequality via the vector infinity norm.
\end{proof}

Lemma \ref{lem:smooth-frame} bounds the Hessian of \(V^\pitheta_f\) by a quantity involving \(\| D_f(\theta) \|_\infty\) and its first and second-order directional derivatives. Since these quantities depend on the \(f\)-function and the policy ratio \(\pitheta(a| s)/\piref(a | s)\), the Hessian is difficult to control uniformly. To address this, in \cref{assump:f-div}, we introduce a general structure that is satisfied by the common \(f\)-divergences we introduce in this paper.
\begin{assumption}\label{assump:f-div}
    Define \(c_\piref:= \min_{s,a}\piref(a| s)\). For the \(f\)-function, there exist constants \( c_{f,1},c_{f,2}>0\), such that for any \(x\in(0,c_{\piref}^{-1})\), it holds \(x\abs{f'(x)}\le c_{f,1}\), and \(x^2\abs{f''(x)}\le c_{f,2}\).
\end{assumption}
    \begin{table*}[t]
        \centering
        \setlength{\tabcolsep}{3pt}
        \caption{Constants of \(f\)-functions that satisfy \cref{assump:f-div}.}
        \begin{tabular}{lcccc}
            \toprule
            \(f\)-divergence      & \(f(x)\)                                                         & \(c_{f,1}\)                                                                 & \(c_{f,2}\)              &\(m_s\)                           \\
            \midrule
            \(\alpha\)-divergence & \(\frac{x^{1-\alpha} - (1-\alpha)x - \alpha}{\alpha(\alpha-1)}\) & \((\alpha c_{\piref})^{-1}\)                                                & \(c_\piref^{\alpha-1}\)  &\((\min_a \piref(a | s))^\alpha\) \\
            reverse KL            & \(x \log x\)                                                     & \(\max\{\frac{1}{e},c_{\piref}^{-1}\log c_{\piref}^{-1}\}+c_{\piref}^{-1}\) & \(c_\piref^{-1}\)        &\(1\)                             \\
            forward KL            & \(- \log x\)                                                     & \(1\)                                                                       & \(1\)                    &\(\min_a \piref(a | s)\)          \\
            JS-divergence         & \(x\log x - (x + 1) \log\frac{x + 1}{2}\)                       & \(\max\{\frac{1}{e},c_{\piref}^{-1}\log2\}\)                                & \(1\)                    &\( \min_a\frac{\piref(a|s)}{1+\piref(a|s)}\)                            \\
            \(\chi^2\)-divergence & \((x - 1)^2\)                                                    & \(2 c_\piref^{-2}\)                                                         & \(2 c_\piref^{-2}\)      &\(\min_a\frac{2}{\piref(a|s)}\)                             \\
            \bottomrule
        \end{tabular}
        \label{tab:f-c-value}
    \end{table*}
    
\begin{remark}\label{rmk:f-c-value}
    The values of  \(c_{f,1}\) and \(c_{f,2}\) for the commonly used \(f\)-divergences are listed in \cref{tab:f-c-value}. Note that this assumption is weaker than directly bounding the derivatives of \(f\), allowing \(f\) to be coercive around 0. The open interval \((0,c_{\piref}^{-1})\) contains the range of \(\frac{\pitheta(a | s)}{\piref(a | s)}\). 
\end{remark}
By utilizing the above assumption with \(x=\frac{\pitheta(a | s)}{\piref(a | s)}\), we can bound the \(D_f(\theta)\) and its first and second-order derivatives. Then we can give the smoothness property of the regularizer as well as the total objective function, which are stated in the lemma below.

\begin{lemma}[Non-uniform Lipschitz smoothness condition]\label{lem:smoothness-pg}
    For any \( \theta, \theta' \in \mathbb{R}^{\dS \dA} \), assume
    \(f\)-function satisfies Assumption \ref{assump:f-div}, then the regularized objective function \(\tilde{V}^{\pitheta}_\lambda(\rho)\) is non-uniform Lipschitz smooth, i.e.,
    \begin{equation*}
    \begin{aligned}
        \left\vert  \tilde{V}_\lambda^{\pi_{\theta'}}\!(\rho) \!-\! \tilde{V}_\lambda^{\pitheta}\!(\rho)\! -\! \langle \nabla_\theta \tilde{V}_\lambda^{\pitheta}(\rho), \theta' \!-\! \theta \rangle \right\vert \! \leq \! \frac{L\!_f(\theta,\theta')}{2} \|\theta' - \theta\|^2,
        \end{aligned}
    \end{equation*}
    where the smoothness factor is a function of \(\theta, \theta'\) given by
    \begin{equation*}
        \begin{aligned}
            L_f(\theta,\theta'):=\frac{8}{(1-\gamma)^3}\Big( & r_{\rm max}+\lambda\sup_{\tilde{\theta}\in [\theta, \theta']}\left\|f\left(w^{\pi_{\tilde{\theta}}}\right)\right\|_{\infty} +\lambda(1-\gamma)\left(c_{f,1}+c_{f,2}/2\right)\Big).
        \end{aligned}
    \end{equation*}
\end{lemma}
\begin{proof}
    We first bound $\left\| D_f(\theta)\right\|_\infty$, $\norm{ \mathcal{D}_\theta D_f(\theta)[v]}_\infty$, and $\norm{\mathcal{D}_\theta^2 D_f(\theta)[v,v]}_\infty$ separately. For any fixed state \(s\), the definition of \(f\)-divergence yields
    \[\abs{D_f(\theta_s)} = \abs{\sum_a \piref(a | s) f \left(w^{\pitheta}_{sa}\right)}\le\max_a \abs{f\left(w^{\pitheta}_{sa}\right)}.\]
    Taking the maximum over all states \(s\in\cS\) gives 
    \begin{equation}
    \norm{D_f(\theta)}_\infty\le \norm{f(w^\pitheta)}_\infty.\label{smo-lem:d-0}
    \end{equation}
    When considering forward KL, the bound is derived as \begin{align*}
        \norm{D_f(\theta)}_\infty &\le \max_s D_f(\theta_s) \\
        & = \max_s\sum_a \piref(a | s)\log\frac{\piref(a | s)}{\pitheta(a | s)} \\
                       & \le \max_s\sum_a \piref(a | s)\abs{\log\pitheta(a | s)} \\
                       & \le \max_{s,a}\abs{\log\pitheta(a | s)} = \|\log\pitheta\|_{\infty}.\numberthis\label{smo-lem:h-0}
    \end{align*}
    To bound the first-order derivative, we compute
    \begin{align*}
           \mathcal{D}_\theta D_f(\theta_s)[v]  & = \sum_a \left(\mathcal{D}_\theta \pi_{\theta}(a | s)[v]\right)f'\left(w^{\pitheta}_{sa}\right) \\
         & = \sum_a  v_s^\top \bH(\theta_s)\delta_a f'\left(w^{\pitheta}_{sa}\right)                                     \\
         & = \sum_a \pitheta(a | s) f'(w^{\pitheta}_{sa}) \Big(v_{sa} - \sum_{a'} v_{sa'}\pitheta(a' | s)\Big) ,
    \end{align*}
    where the second equality uses \eqref{eq:pi-alpha-deriv}. By Assumption \ref{assump:f-div},
    \begin{equation*}
        \abs{f'\left(w^{\pitheta}_{sa}\right)}\le c_{f,1}\frac{\piref(a | s)}{\pitheta(a | s)},
    \end{equation*}
    which implies
    \begin{align*}
         \norm{\mathcal{D}_\theta D_f(\theta)[v]}_\infty = & \max_s \left\vert \mathcal{D}_\theta D_f(\theta_s)[v]\right\vert             \\
         \le & \max_s \sum_a \piref(a | s) c_{f,1} \Big(\abs{v_{sa}} +\sum_{a'} \pitheta(a' | s) \abs{v_{sa'}} \Big) \\
         \le & 2c_{f,1}\max_s \norm{v_s}_\infty \leq 2c_{f,1}\|v\|_2.\numberthis\label{smo-lem:d-1}
    \end{align*}
    For the second-order derivative, it is calculated as
    \begin{align*}
          & \mathcal{D}_\theta^2 D_f(\theta_s)[v,v] \\
          = & \sum_a v_{sa}\left(\mathcal{D}_\theta \pitheta(a|s)[v]\right)\Big(f'\left(w^{\pitheta}_{sa}\right) -\sum_{a'}\pitheta(a' |s) f'\left(w^{\pitheta}_{sa'}\right)\Big)\\
         & + \sum_a v_{sa}\pitheta(a | s)\left(\mathcal{D}_\theta \pitheta(a|s)[v]\right) \frac{f''\left(w^{\pitheta}_{sa}\right)}{\piref(a | s)}                               \\
         & - \sum_a v_{sa}\pitheta(a | s) \sum_{a'}\left(\mathcal{D}_\theta \pitheta(a'|s)[v]\right) f'\left(w^{\pitheta}_{sa'}\right) \\
         & - \sum_a v_{sa}\pitheta(a | s) \sum_{a'}\frac{\pitheta(a' | s)}{\piref(a' | s)} \left(\mathcal{D}_\theta \pitheta(a'|s)[v]\right)f''\left(w^{\pitheta}_{sa'}\right).
    \end{align*}
    Recall that
    \begin{align*}
        \left\vert  \mathcal{D}_\theta \pitheta(a|s)[v] \right\vert & \le \abs{v_{sa}\pitheta(a | s)} + \abs{v_s^\top \pitheta(s)\pitheta(a | s)}\le 2\pitheta(a | s)\|v_s\|_\infty,
    \end{align*}
    with the bound of \(f''\) in \cref{assump:f-div},
    we have
    \begin{align*}
         & \left\vert \mathcal{D}_\theta^2 D_f(\theta_s)[v,v]\right\vert                                                                                                                     \\
         \leq & 2 \sum_a \abs{v_{sa}}\Big\{\norm{v_s}_\infty\Big(c_{f,1}\piref(a| s) + c_{f,1}\pitheta(a | s)  \Big) + \norm{v_s}_\infty c_{f,2} \piref(a | s) \\
         &\qquad +  \pitheta(a | s)c_{f,1} \norm{v_s}_\infty + \pitheta(a | s) c_{f,2} \norm{v_s}_\infty\Big\} \\
         \leq & (6c_{f,1}+ 4c_{f,2}) \norm{v_s}^2_\infty \leq (6c_{f,1}+ 4c_{f,2}) \norm{v_s}^2_2.
    \end{align*}
    Taking the maximum over all states \(s\in\cS\), it follows that
    \begin{equation}
        \norm{\mathcal{D}_\theta^2 D_f(\theta_s)[v,v]}_2\le (6c_{f,1}+ 4c_{f,2})\|v\|_2^2.\numberthis\label{smo-lem:d-2}
    \end{equation}
    Next, using \cref{lem:smooth-frame} together with \eqref{smo-lem:d-0}, \eqref{smo-lem:d-1}, and \eqref{smo-lem:d-2},  for any $\theta,v\in\R^{\dS\dA}$ and $s\in\cS$, we obtain
    \begin{align*}
          & \left\vert  \mathcal{D}_\theta^2 V_f^{\pi_{\theta}}(s) [v,v] \right\vert\\
          \le & \frac{8 \gamma \|v\|_2^2}{(1 - \gamma)^3} \| D_f(\theta) \|_\infty + \frac{4 \gamma \|v\|_2}{(1 - \gamma)^2} \left\| \mathcal{D}_\theta D_f(\theta)[v] \right\|_\infty + \frac{1}{1 - \gamma} \left\|\mathcal{D}_\theta^2 D_f(\theta)[v,v]\right\|_\infty \\
         \le & \frac{\|v\|^2_2}{(1-\gamma)^3}\Big(8\gamma\left\|f\left(w^{\pitheta}\right)\right\|_{\infty}+8\gamma(1-\gamma)c_{f,1} +(1-\gamma)^2(6c_{f,1}+4c_{f,2})\Big) \\
         \le & \frac{8\|v\|^2_2}{(1-\gamma)^3}\left(\left\|f\left(w^{\pitheta}\right)\right\|_{\infty}+(1-\gamma)\left(c_{f,1}+\frac{c_{f,2}}{2}\right)\right).\numberthis\label{eq:smo-vtHv}
    \end{align*}
    For any vectors \( \theta,\theta' \in \mathbb{R}^{\dS \dA} \), and state \( s \in \mathcal{S} \), by Taylor expansion, there exists \( \tilde{\theta} \) on the line segment between \( \theta \) and \( \theta' \) such that
    \begin{align*}
        & \left\vert  V_f^{\pi_{\theta'}}(s) - V_f^{\pitheta}(s) - \langle \nabla_\theta V_f^{\pitheta}(s), \theta' - \theta \rangle \right\vert = \frac{1}{2} \left\vert  \mathcal{D}_\theta^2 V_f^{\pi_{\tilde{\theta}}}(s) [\theta' - \theta,\theta' - \theta] \right\vert\\
       \le & \frac{4 \| \theta' - \theta \|^2_2}{(1 - \gamma)^3}\left(\left\|f\left(w^{\pi_{\tilde{\theta}}}\right)\right\|_{\infty}+(1-\gamma)\left(c_{f,1}+\frac{c_{f,2}}{2}\right)\right)                                             \\
        \le & \frac{4 \| \theta' - \theta \|^2_2}{(1 - \gamma)^3}\!\left(\sup_{\tilde{\theta}\in [\theta, \theta']}\!\left\|f\left(w^{\pi_{\tilde{\theta}}}\right)\right\|_{\infty}+(1-\gamma)\!\left(c_{f,1}+\frac{c_{f,2}}{2}\right)\!\right),
    \end{align*}
    where the first inequality follows from \eqref{eq:smo-vtHv} at point $\tilde{\theta}$ in direction $\theta'-\theta$.
    Combining this with the \( \frac{8r_{\rm max}}{(1 - \gamma)^3} \)-smoothness of \( V^{\pitheta}(\rho) \) completes the proof.
\end{proof}

We analyze the smoothness properties of both forward and reverse KL-regularized value functions as special cases of our general framework. For the forward KL divergence, we obtain a non-uniform smoothness with constants depends only on \(\theta\) and \(\theta'\). For the reverse KL divergence, we generalize the analysis of \cite{mei2020global}, which is only for entropy regularizer, and establish a uniform smoothness condition. The following corollary provides an explicit quantification of the Lipschitz constants.

\begin{corollary}[Smoothness conditions for forward KL and reverse KL]
\label{cor:smoothness-fkl-rkl}
    The forward KL-regularized value function is non-uniform Lipschitz smooth with the smoothness factor \[L_F(\theta,\theta'):= \frac{8 \big( r_{\rm max} + \lambda \max\big\{\|\log\pi_{\theta}\|_\infty,\|\log\pi_{\theta'}\|_\infty\big\}\big)}{(1-\gamma)^3}.
    \]
    The reverse KL-regularized value function \(\tilde{V}^{\pitheta}_\lambda(\rho)\) is uniformly Lipschitz smooth, with smoothness factor \(L_R\). Let \(L_{\rm ref}:=\norm{\log\piref}_\infty\), then
    \[
    L_R:= \frac{8r_{\rm max}}{(1-\gamma)^3}+\lambda\left[\frac{4\gamma(\log\dA+2L_{\rm ref})}{(1-\gamma)^2}+\frac{8\gamma L_{\rm ref}}{(1-\gamma)^3}+\frac{3\log\dA+6L_{\rm ref}+2}{1-\gamma}\right].
    \]
\end{corollary}
\begin{proof}
    For forward KL-regularization, a tighter smoothness factor follows by explicitly evaluating the directional derivatives of \(D_f(\theta)\): 
    \begin{align*}
         \frac{\partial D_f(\theta_s)}{\partial {\theta_{s'a}}}&=  - \sum_{a'}\piref(a' | s)\left[\frac{\partial}{\partial \theta_{s'a}}\log \pi_{\theta}(a' | s)\right] = \mathbbm{1}_{s=s'}\left(\piref(a | s)-\pi_{\theta}(a | s)\right),
    \end{align*}
    Therefore, the following two equations hold:
    \begin{align*}
       & \mathcal{D}_\theta D_f(\theta_s)[v] = \sum_a v_{sa}\frac{\partial D_f(\theta_s)}{\partial {\theta_{s'a}}} = v_s^\top \!\left(\pitheta(s)- \piref(s)\right), \\
        & \mathcal{D}_\theta^2 D_f(\theta_s)[v,v] = D_{\theta_s}^2 D_f(\theta_s)[v_s,v_s] = v_s^\top \mathbf{H}(\theta_s)v_s.
    \end{align*}
    The infinity norms can be further bounded as
    \begin{align}
        & \norm{\mathcal{D}_\theta D_f(\theta)[v]}_\infty =\max_s \abs{\mathcal{D}_\theta D_f(\theta_s)[v]}
        \le \max_s\|v_s\|_\infty\cdot\|\piref(s)-\pitheta(s)\|_1\le 2\|v\|_\infty\le 2\|v\|_2,  \label{smo-lem:h-1} \\
         &\norm{\mathcal{D}_\theta^2 D_f(\theta)[v,v]}_\infty =\max_s \abs{\mathcal{D}_\theta^2 D_f(\theta_s)[v,v]} =\max_s \left\vert  v_s^\top \bH(\theta_s)v_s\right\vert             \le \max_s\|v_s\|_2^2       \le \|v\|_2^2,  \label{smo-lem:h-2}
    \end{align}
    where the first inequality is because the eigenvalues of the positive semi-definite matrix \(\bH(\theta_s)\) are all inside \([0,1]\). Combining the above inequalities \eqref{smo-lem:h-0}, \eqref{smo-lem:h-1}, and \eqref{smo-lem:h-2}, we obtain
    \begin{align*}
          \quad \left\vert  \mathcal{D}_\theta^2 V_f^{\pi_\theta}(s)[v,v]  \right\vert
         & \le \frac{8 \gamma \|v\|_2^2}{(1 - \gamma)^3}  \|\log\pi_{\theta}\|_\infty + \frac{8 \gamma \|v\|_2^2}{(1 - \gamma)^2}+\frac{\|v\|_2^2}{1-\gamma}                                                      \\
         & \le \frac{8 \gamma \|v\|_2^2}{(1 - \gamma)^3} \|\log\pi_{\theta}\|_\infty + \frac{8 \|v\|_2^2}{(1 - \gamma)^2}                                                                                        \le \frac{8 \|v\|_2^2}{(1 - \gamma)^3} \|\log\pi_{\theta}\|_\infty,
    \end{align*}
using \(\min_{s,a} \pitheta(a | s) \leq 1/\abs{\cA}\) and \(\|\log\pi_{\theta}\|_\infty \geq \log \abs{\cA}>1\) when \(\abs{\cA}\gg 1\).     Following the same discussion as in Lemma \ref{lem:smooth-frame}, there exists \(\tilde{\theta}\) on the line segment between \(\theta\) and \(\theta'\) such that
    \begin{align*}
         & \quad \left\vert  \tilde{V}^{\pi_{\theta'}}_\lambda(\rho) - \tilde{V}^{\pitheta}_\lambda(\rho) - \langle \nabla_\theta \tilde{V}^{\pitheta}_\lambda(\rho), \theta' - \theta \rangle \right\vert                  \\
         & \le \frac{4 \left( r_{\rm max} + \lambda \max_{\tilde{\theta}\in[\theta,\theta']}\|\log\pi_{\tilde{\theta}}\|_\infty\right)}{(1 - \gamma)^3} \cdot \| \theta' - \theta \|^2_2 \\
         & \le \frac{4 \left( r_{\rm max} + \lambda \max\big\{\|\log\pi_{\theta}\|_\infty,\|\log\pi_{\theta'}\|_\infty\big\}\right)}{(1 - \gamma)^3} \cdot \| \theta' - \theta \|^2_2,
    \end{align*}
    where the last inequality follows from the monotonicity of the log-function. We conclude the first result for forward KL.

For reverse KL-regularization, we bound the smoothness contribution of \(V_f^\pitheta(\rho)\) directly through \cref{lem:smooth-frame}.
First, from $D_f(\theta_s)\ge 0$ and
\[
    D_f(\theta_s)\le \sum_a\pitheta(a|s)\log\frac{1}{\piref(a|s)}\le \norm{\log\piref}_\infty,
\]
we have
\begin{equation}
\norm{D_f(\theta)}_\infty\le\norm{\log\piref}_\infty\,.\label{smo:rkl-0}
\end{equation}
By the chain rule, we obtain
\begin{align*}
    \frac{\partial D_f(\theta_s)}{\partial\theta_{s',a}}&=\mathbbm{1}_{s=s'}\,\delta_a^\top\mathbf{H}(\theta_s)\big(\log\frac{\pitheta(s)}{\piref(s)}+ \mathbf{1}\big) =\mathbbm{1}_{s=s'}\,\pitheta(a|s)\big(\log\frac{\pitheta(a|s)}{\piref(a|s)}-D_f(\theta_s)\big),
\end{align*}
where in the last equality we use $\mathbf{H}(\theta_s) \mathbf{1}=0$. Denote 
\begin{align*}
    r(\theta_s):=\log\frac{\pitheta(a|s)}{\piref(a|s)}-D_f(\theta_s), \quad  w(\theta_s):=\nabla_{\theta_s}D_f(\theta_s)=\pitheta(s)\odot r(\theta_s),
\end{align*}
where $\odot$ is elementwise multiplication. Then
\begin{align*}
    \norm{\nabla_\theta D_f(\theta_s)}_1 & = \|w(\theta_s)\|_1\\
    &\le \sum_a \pitheta(a|s)\Big(\abs{\log\frac{\pitheta(a|s)}{\piref(a|s)}}+D_f(\theta_s)\Big)\\
    &\le \sum_a \pitheta(a|s)\big(\abs{\log\pitheta(a|s)}+\abs{\log\piref(a|s)}\big)+D_f(\theta_s)\\
    &\le \,\log\dA + 2\norm{\log\piref}_\infty\,.
\end{align*}
Hence, for the directional derivative along an arbitrary $v\in\R^{\dS\dA}$, we have
\begin{align}
    \norm{\mathcal{D}_\theta D_f(\theta)[v]}_\infty = \max_s\abs{\mathcal{D}_\theta D_f(\theta_s)[v]} & \le \max_s \norm{\nabla_\theta D_f(\theta_s)}_1\cdot \|v\|_\infty  \notag \\
    & \le (\log\dA + 2\norm{\log\piref}_\infty)\cdot \|v\|_2. \label{smo:rkl-1}
\end{align}
Using the explicit formula for $\nabla_{\theta_s}^2 D_f(\theta_s)$
\begin{align*}
    [\nabla_{\theta_s}^2 D_f(\theta_s)]_{i,j} &=\frac{\partial w_i(\theta_s)}{\partial \theta_{s,j}}\\
   &=\frac{\partial \pitheta(i|s)}{\partial \theta_{s,j}}r_i(\theta_s) + \pitheta(i|s)\frac{\partial r_i(\theta_s)}{\partial \theta_{s,j}}\\
    &= \pitheta(i|s)\big(\delta_{ij}-\pitheta(j|s)\big)\big(r_i(\theta_s)\!+\!1\big)-\pitheta(i|s)\pitheta(j|s)r_j(\theta_s),
\end{align*}
we obtain
\begin{align}
    & v_s^\top \nabla_{\theta_s}^2 D_f(\theta_s) v_s = \sum_{i,j\in\cA} \left[\nabla_{\theta_s}^2 D_f(\theta_s)\right]_{i,j}\, v_{s,i}v_{s,j}\notag\\
&= \langle w(\theta_s),\, v_s\odot v_s\rangle +\langle \pi_\theta(s),\, v_s\odot v_s\rangle -2\bigl(\pi_\theta(s)^\top v_s\bigr) \big( {w(\theta_s)}^\top v_s\big) -\bigl(\pi_\theta(s)^\top v_s\bigr)^2.\label{smo:rkl-2}
\end{align}
Using the inequality \(\langle v,w\rangle \leq \norm{v}_\infty \norm{w}_1\) together with
\[
\|v_s\odot v_s\|_1 = \|v_s\|_2^2\le\|v\|_2^2 , \quad \|\pi_\theta(s)\|_\infty\le\|\pi_\theta(s)\|_2\le 1,
\]
from \eqref{smo:rkl-2}, we obtain
\begin{align*}
&\bigl|v_s^\top \nabla_{\theta_s}^2 D_f(\theta_s)v_s\bigr|\\
\le& \,(\|w(\theta_s)\|_\infty +\|\pi_\theta(s)\|_\infty)\|v_s\odot v_s\|_1 + 2\|\pitheta(s)\|_2 \|w(\theta_s)\|_2 \,\|v\|_2^2 \, + \|\pitheta(s)\|_2^2 \|v\|_2^2\\
\le& \,(\|w(\theta_s)\|_\infty + 1)\|v\|_2^2 \,+  2\|w(\theta_s)\|_2 \,\|v\|_2^2 \,+ \|v\|_2^2\\
\le & \big(3\|w(\theta_s)\|_1 + 2\big)\,\|v\|_2^2\\
\le & \big(6\norm{\log\piref}_\infty + 3\log\dA + 2\big)\,\|v\|_2^2.
\end{align*}
Consequently,
\begin{equation}
    \norm{\mathcal{D}_\theta^2 D_f(\theta)[v,v]}_\infty\le \big(6\norm{\log\piref}_\infty + 3\log\dA + 2\big)\,\|v\|_2^2.\label{smo:rkl-3}
\end{equation}
Substituting \eqref{smo:rkl-0},\eqref{smo:rkl-1}, and \eqref{smo:rkl-3} into \cref{lem:smooth-frame} gives
\begin{align*}
\left|\mathcal{D}_\theta^2 V_f^{\pitheta}(s)[v,v]\right| &\leq\Bigg( \frac{4\gamma(\log\dA+2\norm{\log\piref}_\infty)}{(1-\gamma)^2} +\frac{8\gamma\norm{\log\piref}_\infty}{(1-\gamma)^3} \\
&\qquad +\frac{3\log\dA+6\norm{\log\piref}_\infty+2}{1-\gamma} \Bigg)\norm{v}_2^2.
\end{align*}
Combining this with the $\frac{8r_{\rm max}}{(1-\gamma)^3}$-smoothness of $V^\pitheta(\rho)$ completes the reverse KL case.
\end{proof}

When \(\piref\) is uniform, reverse KL regularization differs from entropy regularization only by a policy-independent constant. In this case, \(L_{\rm ref}=\log\dA\), and
\(L_R=O\big((r_{\rm max}+\lambda\log\dA)/(1-\gamma)^3\big)\), with the same dependence on \(\dA\) and \(\gamma\) as the entropy-regularized smoothness bound obtained by combining the unregularized bound with Lemma~14 of \cite{mei2020global}. The smoothness constants differ because Corollary~\ref{cor:smoothness-fkl-rkl} is derived through the general \(f\)-divergence framework. The explicit forms above also clarify the respective stepsize restrictions: the forward KL analysis uses the segment-dependent factor \(L_F(\theta,\theta')\), whereas the reverse KL analysis admits the policy-independent factor \(L_R\).

\subsection{Non-uniform Łojasiewicz property of the objective function}
The Łojasiewicz property is the key to guaranteeing the global convergence of PG. As claimed by \cite{mei2020global}, the entropy-regularized value function softmax policy admits a local Łojasiewicz property, with the factor depending on \(\pitheta\). In this section, we will manage to show a similar result for the general \(f\)-divergence regularized objective function.
To derive the desired result, we first need to derive an intermediate lemma that helps connect the gradient \(\nabla_\theta \tilde{V}^{\pitheta}_\lambda(u)\) and the optimality value gap.

\begin{lemma}[Strong concavity of \(F_s^\lambda(p)\) over the simplex]\label{lem:strong-concave-F}
    Define the function \( F_s^\lambda: \Delta_{\abs{\mathcal{A}}} \to \mathbb{R} \) as
\[
    F_s^\lambda(p) := p^\top \tilde{Q}^{\pi_\theta}(s, \cdot) - \lambda D_f(p,\piref(s)),
\]
where \( \tilde{Q}^{\pi_\theta}(s, \cdot) \) is a fixed vector. It is \(\lambda m_s\)-strongly concave on the simplex with 
\[m_s = \inf_{p\in\Delta_{\abs{\cA}}, a \in \mathcal{A}} f'' \left( \frac{p_a}{\piref(a | s)} \right)\frac{1}{\piref(a | s)}\]
and satisfies the gradient domination condition
    \begin{equation*}
        F_s^\lambda(q^*_s) -  F_s^\lambda(\pi(s))\leq \frac{1}{2\lambda m_s}\norm{\Proj_\mathcal{T}(\nabla_{\pi(s)} F_s^\lambda(\pi(s)))}_2^2,
    \end{equation*}
    where \(q^*_s \in \argmax_{q\in \Delta_{\abs{\cA}}} F_s^\lambda(q)\), and \(\Proj_\mathcal{T}(v)\) denotes the projection of \(v\) onto the tangent space of the probability simplex.
\end{lemma}
\begin{proof}
Define \(m_s  := \inf_{p\in\Delta_{\abs{{\cA}}},a \in \mathcal{A}} f''\left(p_a/\piref(a|s) \right) \piref(a | s)^{-1}.\)
    Since \(f\)-function is convex, we know $m_s \geq 0$. By direct computation, we can conclude $F_s^\lambda(p)$ is $\lambda m_s$-strongly concave as 
        \[ \nabla^2_p F_s^\lambda(p) = -\lambda\cdot\diag\left\{\frac{f''\big(p_a/\piref(a|s))}{\piref(a|s)}:a\in\cA\right\}\preceq 0.\]
    Let \( q^*_s \in  \arg \max_{q \in \Delta_{\abs{\mathcal{A}}}} F_s^\lambda(q) \), and define the hyperplane \( \mathcal{T} := \mathbf{1}^\perp \) such  that \(q_s^*+\mathcal{T}\) contains \( \Delta_{\abs{\mathcal{A}}} \) as a subset. Then, for any simplex \(p \in \Delta_{\abs{\cA}}\) it holds that
    \begin{align*}
        F_s^\lambda(q^*_s) - F_s^\lambda(p) \leq & \langle \nabla_{p} F_s^\lambda(p), q^*_s - p \rangle - \frac{\lambda m_s}{2} \|q^*_s - p\|_2^2                          \\
        \leq                                     & \max_{v \in \mathcal{T}} \left\{ \langle \nabla_{p} F_s^\lambda(p), v \rangle - \frac{\lambda m_s}{2} \|v\|_2^2\right\} =                                         \frac{1}{2\lambda m_s} \| \Proj_\mathcal{T} (\nabla_{p} F_s^\lambda(p)) \|_2^2.
    \end{align*}
    Replacing \(p\) with the policy \(\pi(s)\), the final result follows.
\end{proof}

With the strong concavity of the function \(F^\lambda_s(p)\), in the next lemma, we will show the relationship of it to the gradient of the objective function and the optimality gap, then give the local Łojasiewicz property for the \(f\)-divergence regularized value function.

\begin{lemma}[Non-uniform Łojasiewicz]\label{lem:PL-f}
   With Assumptions~\ref{asup:min-piref} and~\ref{asup:min-u}, define $c_m:=\min_s m_s$, the regularized value function satisfies 
    \begin{equation*}
    \begin{aligned}
        & \| \nabla_\theta \tilde{V}^{\pitheta}_\lambda(u) \|_2^2   \geq 2\lambda c_u c_m\big( \min_{a, s} \pitheta(a | s) \big)^2      \big\| d^{\pi^*}_\rho / d^{\pitheta}_u \big\|_{\infty}^{-1} \big( \tilde{V}^{\pi^*}_\lambda(\rho) - \tilde{V}^{\pitheta}_\lambda(\rho) \big).
        \end{aligned}
    \end{equation*}
    The strong concavity coefficient \(m_s\) for the commonly used \(f\)-divergence is provided in \cref{tab:f-c-value}.
\end{lemma}
\begin{proof}
     Rewrite the Performance Difference Lemma \ref{lem:performance-difference} with the function \(F_s^\lambda\) defined in \cref{lem:strong-concave-F}, we can obtain
    \begin{align*}
         \tilde{V}^{\pi^*}_\lambda(\rho) - \tilde{V}^{\pi_\theta}_\lambda(\rho)          & = \frac{1}{1-\gamma} \sum_s d^{\pi^*}_{\rho}(s) \left( F_s^\lambda(\pi^*(s)) - F_s^\lambda(\pitheta(s)) \right) \\
         & \leq  \frac{1}{1-\gamma} \sum_s d^{\pi^*}_{\rho}(s) \left( F_s^\lambda(q^*_s) - F_s^\lambda(\pitheta(s)) \right)        \\
         & \le \frac{1}{2\lambda \min_s m_s(1-\gamma)} \sum_s d^{\pi^*}_{\rho}(s) \|Z_\theta(s) \|_2^2 \\
        & \leq  \frac{1}{2\lambda c_m(1-\gamma)} \norm{\frac{d^{\pi^*}_\rho}{d^{\pitheta}_u}}_\infty \sum_s d^{\pitheta}_u(s)\norm{Z_\theta(s)}_2^2 ,\numberthis \label{eq:gap-z-ineq}
    \end{align*}
    where \(Z_\theta(s) := \Proj_{\cT} (\nabla_{\pi_\theta(s)} F_s^\lambda (\pitheta(s)))\) denotes the projected vector, and it takes the form
        $$Z_\theta(s) = \left(\mathbf{I} - \frac{1}{\abs{\cA}}\mathbf{1}\mathbf{1}^\top\right) \left(\tilde Q^{\pitheta}(s,\cdot) - \lambda f'\left( \frac{\pitheta(s)}{\piref(s)} \right)\right) .$$
    By Lemma 23 in \cite{mei2020global}, for any \(v \in \cT\), it holds \(\norm{\bH(\theta_s) v}_2 \geq \min_a \pitheta(a | s) \|v\|_2\). With the fact that \(\mathbf{H}(\pi)\mathbf{1} =0\),  we can lower bound the gradient as
    \begin{align*}
        \left\|\frac{\partial \tilde V^{\pitheta}_\lambda(u)}{\partial \theta_s}\right\|_2 =  \frac{1}{1-\gamma} d^{\pitheta}_u(s) \left\| \bH(\theta_s) Z_\theta(s)\right\|_2 \geq                                                                                 \frac{1}{1-\gamma}d^{\pitheta}_u(s)  \min_a \pitheta(a | s)\left\|  Z_\theta(s)\right\|_2.
    \end{align*}
    Summing over all states, it holds
\begin{equation}\label{eq:grad-Z-ineq}
    \begin{aligned}
        & \left\|\nabla_\theta \tilde V^{\pitheta}_\lambda(u) \right \|_2^2 \geq \frac{1}{(1-\gamma)^2} \sum_s d^{\pitheta}_u(s)^2 (\min_{a}\pitheta(a | s))^2 \left\|Z_\theta(s) \right\|_2^2.
        \end{aligned}
    \end{equation}
    Combining inequalities \eqref{eq:gap-z-ineq} and \eqref{eq:grad-Z-ineq}, we can show the claim through the following steps.
    \begin{align*}
            \norm{\nabla_\theta \tilde V^{\pitheta}_\lambda(u)}_2^2   & \geq  \frac{1}{(1-\gamma)^2} \sum_s d^{\pitheta}_u(s)^2 (\min_{a}\pitheta(a | s))^2 \left\|Z_\theta(s) \right\|_2^2                                                                                                          \\
        & \geq  \frac{\min_s d^{\pitheta}_u(s)}{(1-\gamma)^2} (\min_{a,s}\pitheta(a | s))^2 \sum_s d^{\pitheta}_u (s)\norm{Z_\theta(s)}_2^2                                                                                           \\
        & \geq  \frac{2\lambda }{1-\gamma}   c_m   \min_s d^\pitheta_u(s) (\min_{s,a}\pitheta(a | s))^2 \norm{d^{\pi^*}_\rho/d^\pitheta_u}^{-1}_\infty (\tilde V^{\pi^*}_\lambda(\rho) - \tilde V^\pitheta_\lambda (\rho)) \\
       &  \geq  2\lambda c_u  c_m (\min_{s,a}\pitheta(a | s))^2  \norm{d^{\pi^*}_\rho / d^\pitheta_u }^{-1}_\infty (\tilde V^{\pi^*}_\lambda(\rho) - \tilde V^\pitheta_\lambda (\rho)). \qedhere
    \end{align*}
\end{proof}

\section{Theoretical guarantees of PPO-Clip for KL-regularized function}
Equipped with the aforementioned properties, we can provide the theoretical guarantees of the PPO-Clip algorithm. The geometries of \(f\)-divergences vary significantly, leading to diverse analysis techniques. Therefore, in this section, we will focus specifically on the KL-divergence, including both forward and reverse KL.
Whenever \(\pi_{n,k}=\pi_{\theta_{n,k}}\), we use the convention
$\nabla_\theta \tilde V_\lambda^{\pi_{n,k}}(u)
:=\left.\nabla_\theta \tilde V_\lambda^{\pi_\theta}(u)\right|_{\theta=\theta_{n,k}}$.
For later use, we also define the clipping factor
\(\kappa_\varepsilon := 1/\log(1+\varepsilon_h) -1/\log(1-\varepsilon_l).\)

\subsection{PPO-Clip for forward KL}

The surrogate objective function with the forward KL-regularizer is defined as
\begin{align*}
     \mathcal{L}_{n}(\theta) &:= \E_{\pi_{n,1}}  \Big[\min\big\{r_\theta^n(s,a) \tilde{A}^{\pi_{n,1}}_\lambda(s,a), \text{clip}\big(r_\theta^n(s,a);1-\varepsilon_l,  1+\varepsilon_h\big) \tilde{A}^{\pi_{n,1}}_\lambda(s,a)\big\} \\
&\qquad +  \lambda  r_{\rm ref}^n(s,a) \log \pi_{\theta}(a | s) \Big].
\end{align*}
When analyzing the theoretical property of the RL algorithms, a common assumption on the bounded reward \(r(s,a)\) implies the bounded advantage function. However, for the forward KL-divergence regularized objective function considered in this section, the corresponding reward and the advantage function are not bounded in general.  In the following lemma, we give a scenario in which the regularized advantage function is guaranteed to be bounded.
\begin{lemma}[Bounded advantage function]\label{lem:bounded-advantage} With Assumptions~\ref{asup:min-piref} and~\ref{asup:min-u}, for an arbitrary policy \(\pitheta\) such that \(\tilde{V}^\pitheta_\lambda(u) \geq \tilde{V}_0\), where \(\tilde{V}_0\) is some predetermined constant, then there exist constants \(c_{\pitheta}^F>0\), \(C_A<\infty\) such that  \(\min_{a,s}\pitheta(a | s) \geq c_\pitheta^F\), and \(\tilde{A}_{\rm max}^F:= \sup_{s,a}|\tilde{A}^\pi(s,a)|\leq \frac{1}{1-\gamma}(r_{\rm max} + \lambda  C_A) \).
\end{lemma}
\begin{proof}
By the definition, \(\tilde{V}^\pitheta_\lambda(u) \geq \tilde{V}_0\) implies
    \begin{align*}
    \lambda     V^\pitheta_f(u) \leq  \E_{u,\tau}\Big[\sum_{t=0}^\infty \gamma^t r(s_t,a_t) \Big] - \tilde{V}_0 \leq                 \frac{r_{\rm max}}{1-\gamma} - \tilde{V}_0.
    \end{align*}
    Also, since \(V^\pitheta_f(u) = \frac{1}{1-\gamma}\sum_s d^\pitheta_u(s) D_{KL}(\piref(s) , \pitheta(s))\) and \(d^\pitheta_u(s)  \geq (1-\gamma) c_u\),
 we can derive
    \begin{align*}
        \frac{r_{\rm max}}{1-\gamma} - \tilde{V}_0 \geq  \lambda  \sum_s u(s)D_{KL}(\piref(s) , \pitheta(s))
        \geq                                                                           \lambda c_u \sup_s D_{KL}(\piref(s) , \pitheta(s)).
    \end{align*}
    For any \(s\), it holds \[D_{KL}(\piref(s) , \pitheta(s)) \leq C_A:= \frac{1}{\lambda c_u}\left(\frac{r_{\rm max}}{1-\gamma}-\tilde{V}_0 \right), \] and
    \begin{align*}
             D_{KL}(\piref(s) ,\pitheta(s)) &=     \sum_a \piref(a | s) \log \piref(a | s) - \sum_a \piref(a | s) \log \pitheta(a | s)
        \\
        &\geq  - \log\abs{\cA} - c_{\piref}\log  \min_a \pitheta(a | s).
    \end{align*}
    This implies \(\pitheta(a | s)\) is lower bounded as
    \begin{align*}
        \min_a \pitheta(a | s) \geq \exp \left( (-C_A -\log \abs{\cA}) / c_{\piref} \right), \quad \forall s \in \cS.
    \end{align*}
    To bound the advantage \(A^\pitheta_f(s,a)\), by its definition, for any state \(s\) and action \(a\), we have
    \begin{align*}
        A^\pitheta_f(s,a) &=  Q^\pitheta_f(s,a) - V^\pitheta_f(s)                                                                      \\
        &=                    \gamma \sum_{s'}\left(\cP(s' | s,a) - \sum_{a'}\pitheta(a' | s)\cP(s' | s,a') \right)V^\pitheta_f(s')    \\
       &  =                    \gamma \langle P(\cdot | a,s) - \sum_{a'}\pitheta(a' | s) \cP(\cdot | s,a'), V^\pitheta_f(\cdot)\rangle.
    \end{align*}
    Taking the absolute value on both sides, we have \(\abs{A^\pitheta_f(s,a) }\leq  \gamma \norm{V^\pitheta_f(\cdot)}_{\rm span} \leq  \frac{\gamma}{1-\gamma}C_A\),    where the first inequality holds because for two simplexes \(p\), \(q\) and vector \(v\), we can show 
    \begin{align*}
        \abs{\langle p-q,v\rangle}\leq \frac{1}{2}\norm{p-q}_1 \|v\|_{\rm span} \leq \|v\|_{\rm span}.
    \end{align*}
    The last inequality holds because \( D_{KL}(\piref(s),\pitheta(s)) \leq C_A\) and \(V^\pi_f(s) \leq \sum_{t=0}^\infty \gamma^t C_A = \frac{C_A}{1-\gamma}\).
With \(\tilde{A}_{\rm max}^F \leq \sup_{s,a}\abs{A^{\pi}(s,a)}+\lambda \sup_{s,a}\abs{A^{\pi}_f(s,a)}\), we can derive the bound of the advantage function.
\end{proof}

Next, we give a bound of the decreased objective value when one outer-loop of PPO-Clip is applied. 

\begin{lemma}[Descent lemma under forward KL]\label{lem:descent-lemma}For softmax policy parameterization, if \(\tilde{V}^{\pi_{n,1}}_\lambda(u) \geq \tilde{V}_0\), the stepsize satisfies \(S_{\max}:=\sum_{i=1}^{K}\eta_{n,i} \leq \min\{\frac{1-\gamma}{8C_e},  \frac{(1-\gamma) \log 2}{4\tilde{A}_{\rm max}^F+8\lambda}\} \), where  \(C_e:=\tilde{A}_{\rm max}^F(\sqrt{2}\kappa_\varepsilon+3)+\lambda \), then the iteration sequence \(\{\theta_{n,k}\}\) generated by \cref{alg:PPO-Clip} satisfies
    \begin{align*}
              & (\tilde{V}^{\pi^*}_\lambda(u) - \tilde{V}^{\pitheta_{n,K+1}}_\lambda(u)) - (\tilde{V}^{\pi^*}_\lambda(u) - \tilde{V}^{\pi_{n,1}}_\lambda(u)) \\
        &\quad \leq -\frac{3}{4} S_{\max} \|\nabla_\theta \tilde{V}^{\pi_{n,1}}_\lambda(u)\|^2_2 + 2L S_{\max}^2 \|\nabla_\theta \tilde{V}^{\pi_{n,1}}_\lambda(u)\|^2_2 ,
    \end{align*}
    with \(L=\frac{8\left(r_{\rm max}+\lambda\log2 +\lambda\frac{\log\dA+C_A}{\min_{sa}\piref(a | s)}\right)}{(1-\gamma)^3}\).
\end{lemma}

\begin{proof}
    Let \(r_{n,k}(s,a):= \pitheta_{n,k}(a| s)/ \pi_{n,1}(a| s)\), the non-clip region \(\mathcal{N}_{n,k}\) and the clip region \(\mathcal{C}_{n,k}\) are defined as
    \begin{align*}
         & \mathcal{N}_{n, k}(s,a)  :=    \left\{\tilde{A}^{\pi_{n,1}}_\lambda(s,a) \geq 0, r_{n, k}(s,a) \leq 1+\varepsilon_h\right\} \cup \left\{\tilde{A}^{\pi_{n,1}}_\lambda (s,a)<0, r_{n, k}(s,a) \geq 1-\varepsilon_{\ell}\right\}, \\
         & \mathcal{{C}}_{n,k}(s,a) :=  \mathcal{N}_{n, k}^c(s,a).
    \end{align*}
  Then the gradient of the PPO-Clip subproblem \(g_{n,k}\) is 
    \begin{align*}
        g_{n,k} = & \E_{\pi_{n,1}}\Big[  \psi_{n,k}(s,a) \Big(\mathbbm{1}_{\mathcal{N}_{n,k}}(s,a) r_{n,k}(s,a) \tilde{A}^{\pi_{n,1}}_\lambda(s,a) + \lambda r_{\rm ref}^n(s,a) \Big)\Big].
    \end{align*}
    From \cref{lem:gradient}, the gradient of \(\tilde{V}^{\pi_{n,1}}_\lambda(u)\) is
    \begin{align*}
           \nabla_\theta \tilde V^{\pi_{n,1}}_\lambda (u)
        =  \E_{\pi_{n,1}}\Big[\psi_{n,1}(s,a) \big(\tilde{A}^{\pi_{n,1}}_\lambda(s,a) + \lambda r_{\rm ref}^n(s,a)\big)\Big].
    \end{align*}
The PPO-Clip gradient ascent direction \(g_{n,k}\) can be treated as an inexact gradient as follows:
    \begin{align*}
        g_{n,k} =&  \nabla_\theta \tilde V^{\pi_{n,1}}_\lambda(u)  \underbrace{-\E_{\pi_{n,1}}\!\big[ \mathbbm{1}_{\mathcal{C}_{n,k}}(s,a) \,\tilde{A}^{\pi_{n,1}}_\lambda(s,a)\,\psi_{{n,1}}(s,a)\big]}_{:=X_{n,k}}                   \\
                  & + \underbrace{\E_{\pi_{n,1}}\big[\mathbbm{1}_{\mathcal{N}_{n,k}}(s,a)\frac{\nabla_\theta \pitheta_{n,k}(a | s) - \nabla_\theta \pi_{n,1}(a | s)}{\pi_{n,1}(a | s)}\tilde{A}^{\pi_{n,1}}_\lambda(s,a) \big]}_{:=Y_{n,k}} \\
                  & + \underbrace{\E_{\pi_{n,1}}\big[\lambda r_{\rm ref}^n(s,a)(\psi_{{n,k}} (s,a)-\psi_{{n,1}} (s,a))\big]}_{:=Z_{n,k}}.
    \end{align*}
    Next, we need to bound three terms \(X_{n,k}, Y_{n,k}\), and \(Z_{n,k}\) in the above inequality individually. For the first term, by the property of the indicator function, we have
    \begin{align*}
       \mathbbm{1}_{A_{n,1}\geq 0, r_{n,k}\geq 1 + \varepsilon_h} &=    \mathbbm{1}_{A_{n,1}\geq 0, \log r_{n,k}\geq \log( 1 + \varepsilon_h)} \leq                                                             \mathbbm{1}_{ \log r_{n,k}\geq \log( 1 + \varepsilon_h)} \\
       &\leq                                                             \frac{\max \{0,\log r_{n,k}(s,a)\}}{\log (1+\varepsilon_h)} ,            \\
       \mathbbm{1}_{A_{n,1}  \leq 0, r_{n,k}\leq 1  - \varepsilon_l} &=  \mathbbm{1}_{A_{n,1} \leq 0, \log r_{n,k}\leq \log (1- \varepsilon_l) } \leq                                                             \mathbbm{1}_{\log r_{n,k} \leq \log(1 - \varepsilon_l)} \\
       &\leq   \frac{\max\{0, -\log r_{n,k}(s,a)\}}{-\log(1-\varepsilon_l)}.
    \end{align*}
    With these two inequalities, we can derive
    \(\mathbbm{1}_{\mathcal{C}_{n,k}}(s,a) \leq \kappa_\varepsilon\abs{\log r_{n,k}(s,a)}.\)
    Define \(\Delta_{n,k}:= \theta_{n,k} - \theta_{n,1} \), \(\theta(l) :=\theta_{n,1} + l\Delta_{n,k} \), by the Fundamental Theorem of Calculus, it holds
    \begin{align*}
    \abs{\log r_{n,k}(s,a)} =
        \abs{\int_0^1 \psi_{\theta(l)} (s,a)^\top \Delta_{n,k} dl} \leq \norm{\Delta_{n,k}}_2 \int_0^1 \norm{\psi_{\theta(l)} (s,a)}_2 dl \leq \sqrt{2}\norm{\Delta_{n,k}}_2.
    \end{align*}
    Hence, \(\norm{X_{n,k}}_2\) can be bounded as
    \begin{align*}
        \norm{X_{n,k}}_2 &\leq  \E_{\pi_{n,1}}\Big[\mathbbm{1}_{\mathcal{C}_{n,k}}(s,a)\norm{\psi_{n,1} (s,a)} _2\abs{\tilde{A}^{\pi_{n,1}}_\lambda(s,a)} \Big]  \leq                  \frac{\sqrt{2}\tilde{A}_{\rm max}^F}{1-\gamma}\kappa_\varepsilon\norm{\Delta_{n,k}}_2.
    \end{align*}
    To bound \(\norm{Y_{n,k}}_2\), with \eqref{ineq:regularity-pitheta-smooth} in \cref{lem:bounded-psi}, we have
    \begin{align*}
        & \norm{Y_{n,k}}_2
         \leq                \E_{\pi_{n,1}}\Big[1/ {\pi_{n,1}(a | s)} |\tilde{A}^{\pi_{n,1}}_\lambda(s,a)| \big\| \nabla_\theta \pitheta_{n,k}(a | s)    - \nabla_\theta \pi_{n,1}(a | s)\big\| _2 \Big]                                  \leq                \frac{3}{1-\gamma} \tilde{A}_{\rm max}^F\norm{\Delta_{n,k}}_2.
    \end{align*}
    For the last term, with \cref{lem:bounded-psi}, \(Z_{n,k}\) can be bounded as
    \begin{align*}
             \norm{Z_{n,k}}_2
        \leq&  \frac{\lambda}{1-\gamma} \! \E_{s\sim d^{\pi_{n,1}}_u}\! \Big[\sum_a \piref(a | s) \norm{\psi_{\theta_{n,k}} (s,a) - \psi_{n,1} (s,a)}_2 \Big] \\
        \leq&  \frac{\lambda}{1-\gamma}\! \E_{s\sim d^{\pi_{n,1}}_u} \! \Big[\sum_a \piref(a | s) \norm{\Delta_{n,k}}_2 \Big] =    \frac{\lambda}{1-\gamma}\norm{\Delta_{n,k}}_2.
    \end{align*}
    Combining all these bounds, the  direction $g_{n,i}$ is bounded as
    \begin{align*}
        \norm{g_{n,i}}_2 \leq \| \nabla_\theta \tilde{V}^{\pi_{n,1}}_\lambda(u)\|_2 +\frac{C_e}{1-\gamma}\norm{\Delta_{n,i}}_2.
    \end{align*}
    By the definition of \(\norm{\Delta_{n,k}}_2\), and the update rule of the PPO-Clip algorithm, for any \(1\leq k\leq K+1\), it holds
    \begin{align*}
         \norm{\Delta_{n,k}}_2 = \norm{\theta_{n,k}-\theta_{n,1}}_2 = \norm{\sum_{i=1}^{k-1}(\theta_{n,i+1}-\theta_{n,i})}_2 
         \leq \sum_{i=1}^{k-1}\eta_{n,i}\norm{g_{n,i}}_2 \leq   \max_{i\leq k-1} \norm{g_{n,i}}_2\sum_{i=1}^{k-1}\eta_{n,i}.
    \end{align*}
    Take the maximum of \(\norm{g_{n,i}}_2\) among all inner iterations,
\begin{align*}
        \max_{i\leq K} \norm{g_{n,i}}_2 &\leq  \|\nabla_\theta \tilde{V}^{\pi_{n,1}}_\lambda(u)\|_2 +\frac{C_e}{1-\gamma}  \max_{i\leq K} \norm{\Delta_{n,i}}_2                    \\
       &\leq  \|\nabla_\theta \tilde{V}^{\pi_{n,1}}_\lambda(u)\|_2 + \frac{C_e}{1-\gamma} \max_{i\leq K} \norm{g_{n,i}} \sum_{i=1}^{K}\eta_{n,i}.
    \end{align*}
    From the above inequality for \(\max_{i\leq K} \norm{g_{n,i}}_2\), we have
    \begin{align*}
        \Big(1- \frac{C_e}{1-\gamma}  \sum_{i=1}^{K}\eta_{n,i} \Big)\max_{i\leq K}\norm{g_{n,i}}_2 \leq \|\nabla_\theta \tilde{V}^{\pi_{n,1}}_\lambda(u)\|_2.
    \end{align*}
    With the assumption \(\frac{C_e}{1-\gamma} \sum_{i=1}^{K}\eta_{n,i} \leq 1/2\), the ascent direction can be bounded as
    \begin{equation*}
        \max_{i\leq K}\norm{g_{n,i}}_2 \leq 2 \|\nabla_\theta \tilde{V}^{\pi_{n,1}}_\lambda(u)\|_2.
    \end{equation*}
    Also, for the softmax policy, the gradient is bounded as
    \begin{align*}
             \|\nabla_\theta \tilde{V}^{\pi_{n,1}}_\lambda(u)\|_2 \leq& \E_{\pi_{n,1}}\Big[\norm{\psi_{n,1}(s,a)}_2 \abs{\tilde{A}^{\pi_{n,1}}(s,a)} \Big] \\
        &+  \frac{\lambda}{1-\gamma} \E_{s \sim d^{\pi_{n,1}}_u}\Big[ \sum_a \norm{\psi_{n,1}(s,a)}_2  \abs{\piref(a | s) -\pi_{n,1}(a | s)} \Big] \\
        \leq&  \frac{1}{1-\gamma} \sqrt{2}\tilde{A}_{\rm max}^F+ \frac{\sqrt{2} \lambda}{1-\gamma}\E_{s \sim d^{\pi_{n,1}}_u} [\norm{\piref(s) - \pi_{n,1}(s)}_1]                                                  \\
        \leq&  \frac{\sqrt{2}}{1-\gamma} \tilde{A}_{\rm max}^F+ \frac{2\sqrt{2}\lambda}{1-\gamma}.
    \end{align*}
    Recall that \(\abs{\log r_{n,k}(s,a)}\leq \sqrt{2} \norm{\Delta_{n,k}}_2\), hence by setting \(S_{\max}:=\sum_{i=1}^{K}\eta_{n,i} \leq \frac{(1-\gamma) \log 2}{4\tilde{A}_{\rm max}^F+8\lambda}\)    for any \(s,a\) and \(1\leq k\leq K+1\) we can guarantee
    \begin{align*}
        \pitheta_{n,k}(a | s)\geq \frac{1}{2}\pi_{n,1}(a | s) \geq \frac{1}{2}\exp\left(\frac{-\log \abs{\cA} - C_A}{c_{\piref}} \right).
    \end{align*}
As \(\pitheta_{n,k}(a | s)\) is lower bounded in the inner loop, the parameterized value function \(\tilde{V}_\lambda^{\pi_\theta}(u)\) is \(L\)-smooth in the minimum region that contains the trajectory of the inner loop, with \(L=\frac{8\left(r_{\rm max}+\lambda\log2 +\lambda\frac{\log\dA+C_A}{c_{\piref}}\right)}{(1-\gamma)^3}\). Applying this smoothness bound to the outer update gives
    \begin{align*}
             & \left(\tilde{V}_\lambda^{\pi^*}(u)-\tilde{V}_\lambda^{\pi_{n+1,1}}(u)\right) -\left(\tilde{V}_\lambda^{\pi^*}(u)-\tilde{V}_\lambda^{\pi_{n,1}}(u)\right)\\
        \leq & -\left\langle \nabla_\theta \tilde{V}_\lambda^{\pi_{n,1}}(u), \sum_{i=1}^{K}\eta_{n,i}g_{n,i}\right\rangle +\frac{L}{2}\norm{\Delta_{n,K+1}}_2^2\\
        =    & -S_{\max}\|\nabla_\theta \tilde{V}_\lambda^{\pi_{n,1}}(u)\|^2_2 +\frac{L}{2}\norm{\Delta_{n,K+1}}_2^2 -\sum_{i=1}^{K}\eta_{n,i}\left\langle \nabla_\theta \tilde{V}_\lambda^{\pi_{n,1}}(u),X_{n,i}+Y_{n,i}+Z_{n,i}\right\rangle\\
        \leq & -S_{\max}\|\nabla_\theta \tilde{V}_\lambda^{\pi_{n,1}}(u)\|^2_2 +\frac{L}{2}\norm{\Delta_{n,K+1}}_2^2 +\frac{C_e}{1-\gamma}\|\nabla_\theta \tilde{V}_\lambda^{\pi_{n,1}}(u)\|_2 \sum_{i=1}^{K}\eta_{n,i}\norm{\Delta_{n,i}}_2\\
        \leq & -S_{\max}\|\nabla_\theta \tilde{V}_\lambda^{\pi_{n,1}}(u)\|^2_2 +2L S_{\max}^2\|\nabla_\theta \tilde{V}_\lambda^{\pi_{n,1}}(u)\|^2_2 +\frac{2C_e}{1-\gamma}S_{\max}^2 \|\nabla_\theta \tilde{V}_\lambda^{\pi_{n,1}}(u)\|^2_2 .
    \end{align*}
    Under the condition \(S_{\max}\leq (1-\gamma)/(8C_e)\), the proof is complete.
\end{proof}

\cref{lem:descent-lemma} establishes key properties for the inner-loop of PPO-Clip (\cref{alg:PPO-Clip}). Given a suitable initial policy and stepsize, it is guaranteed that the generated policies are lower-bounded and that the regularized value function is locally Lipschitz smooth. Equipped with this lemma, we are ready to show the property of the algorithm for the whole trajectory. First, we will give the stationary convergence in \cref{thm:stationary-converge}, then we will show the global convergence of the algorithm in \cref{thm:global-converge}.

\begin{theorem}[Stationary convergence rate of forward KL]\label{thm:stationary-converge}
    Assume that \cref{alg:PPO-Clip} starts from \(\pitheta_{1,1}=\piref\), and that the stepsize satisfies \(S_{\max}:=\sum_{i=1}^{K}\eta_{n,i} \leq \min\{\frac{1-\gamma}{8C_e},  \frac{(1-\gamma) \log 2}{4\tilde{A}_{\rm max}^F+8\lambda}, \frac{1}{8L}\}  \).  Let \(\tilde{V}_0(u):= \tilde{V}^{\pitheta_{1,1}}_\lambda(u)\), then
    \begin{align*}
        \min_{1\leq n\leq N}\norm{\nabla_\theta \tilde{V}_\lambda^{\pi_{n,1}}(u)}^2_2 \leq \frac{2(\tilde{V}^{\pi^*}_\lambda(u) -\tilde{V}^{\pitheta_{1,1}}_\lambda(u)) }{N S_{\max}}.
    \end{align*}
\end{theorem}
\begin{proof}
    At the first iteration, \(\tilde{V}^{\pi_{n,1}}_\lambda(u)\geq \tilde{V}_0\) holds, and thus \cref{lem:bounded-advantage} is applicable. With the stepsize satisfying \(S_{\max}\leq \min\{\frac{1-\gamma}{8C_e},  \frac{(1-\gamma) \log 2}{4\tilde{A}_{\rm max}^F+8\lambda}, \frac{1}{8L}\}\), \cref{lem:descent-lemma} guarantees \(\tilde{V}^{\pitheta_{n+1,1}}_\lambda(u)\geq \tilde{V}^{\pitheta_{n,1}}_\lambda(u)\geq \cdots \tilde{V}^{\pitheta_{1,1}}_\lambda(u)= \tilde{V}_0(u)\) and
    \begin{align*}
        &(\tilde{V}^{\pi^*}_\lambda(u) - \tilde{V}^{\pitheta_{n+1,1}}_\lambda(u)) - (\tilde{V}^{\pi^*}_\lambda(u) - \tilde{V}^{\pi_{n,1}}_\lambda(u)) \leq  -\frac{1}{2}S_{\max}\|\nabla_\theta \tilde{V}^{\pi_{n,1}}_\lambda(u)\|^2_2\, , \quad \forall n\leq N.
    \end{align*}
    Telescoping the above inequality from \(n=1\) to \(N+1\) gives
    \begin{align*}
        N \min_{1\leq n\leq N} \|\nabla_\theta \tilde{V}^{\pi_{n,1}}_\lambda(u)\|^2_2 \leq \sum_{n=1}^{N}\|\nabla_\theta \tilde{V}^{\pi_{n,1}}_\lambda(u)\|^2_2 \leq \frac{2(\tilde{V}^{\pi^*}_\lambda(u) - \tilde{V}^{\pitheta_{1,1}}_\lambda(u))}{S_{\max}}.
    \end{align*}
    Dividing \(N\) on both sides proves the final result.
\end{proof}

The trajectory-wise policy lower bound makes the non-uniform Łojasiewicz inequality uniform along the iterates. Combined with the descent estimate, it gives a geometric contraction of the value gap. Thus, \cref{thm:global-converge} below establishes global linear convergence for the true initial distribution \(\rho\), while the gradient is evaluated using \(u\).

\begin{theorem}[Global  convergence rate under forward KL]\label{thm:global-converge}
    Assume that \cref{alg:PPO-Clip} starts from \(\pitheta_{1,1}=\piref\), and that the stepsize satisfies \(S_{\max}:=\sum_{i=1}^{K}\eta_{n,i} \leq \min\{\frac{1-\gamma}{8C_e},  \frac{(1-\gamma) \log 2}{4\tilde{A}_{\rm max}^F+8\lambda}, \frac{1}{8L}\}  \). Then, for every \(n\geq1\),
    \begin{align*}
        \tilde{V}^{\pi^*}_\lambda(\rho) - \tilde{V}^{\pi_{n,1}}_\lambda(\rho) \leq \frac{\tilde{V}^{\pi^*}_\lambda(u) - \tilde{V}^{\pitheta_{1,1}}_\lambda(u)}{c_u(1-\gamma)} {e^{-C(n-1)}},
    \end{align*}
    with \(C:=(1-\gamma)\lambda S_{\max} c_u c_{\piref}(c_{\pitheta}^F)^2\norm{d^{\pi^*}_u/u}^{-1}_\infty\).
\end{theorem}

\begin{proof}
    By \cref{cor:value-sub-optimality}, replace the \(f\)-divergence with forward KL-divergence, and define the terms inside the bracket as
    \begin{align*}
        W(s):=&\sum_a (\pi^*(a | s) - \pi_{n,1}(a | s)) \tilde{Q}^{\pi^*}_\lambda(s,a) + \lambda\sum_a \piref(a | s) \log \frac{\pi^*(a | s)}{\pi_{n,1}(a | s)}.
    \end{align*}
    Hence,
    \begin{align*}
       \tilde{V}^*_\lambda(\rho) - \tilde{V}^{\pi_{n,1}}_\lambda(\rho) &= \frac{1}{1-\gamma}\sum_s d^{\pi_{n,1}}_\rho(s) W(s)                                                                          \\
       & =                                                                       \frac{1}{1-\gamma}\sum_s \frac{d^{\pi_{n,1}}_\rho(s)}{d^{\pi_{n,1}}_u(s)} d^{\pi_{n,1}}_u(s) W(s)                  \\
        &\leq                                                                    \frac{1}{c_u(1-\gamma)^2}  \sum_s  d^{\pi_{n,1}}_u(s) W(s)                                             \\
        & =                                                                       \frac{1}{c_u(1-\gamma)}  \left( \tilde{V}^{\pi^*}_\lambda(u) - \tilde{V}^{\pi_{n,1}}_\lambda(u) \right) \numberthis \label{eq:change-initial-distribution},
    \end{align*}
    where the first inequality holds because \(d^{\pi_{n,1}}_\rho(s)\leq 1\), \(d^{\pi_{n,1}}_u(s)\geq (1-\gamma) u(s)\), and \(W(s)\geq0\) for all states \(s\).
    Denote \(\delta_n:=\tilde{V}^{\pi^*}_\lambda(u) - \tilde{V}^{\pi_{n,1}}_\lambda(u)  \), according to \cref{lem:descent-lemma,lem:PL-f}, with \(S_{\max}\leq \min\{\frac{1-\gamma}{8C_e},\frac{(1-\gamma)\log 2}{4\tilde{A}_{\rm max}^F+8\lambda}, \frac{1}{8L}\}\), it holds
    \begin{align*}
         \delta_{n+1}-\delta_n         \leq & -\frac{1}{2}S_{\max}\norm{\nabla_\theta\tilde{V}^{\pi_{n,1}}_\lambda(u)}^2                                                                                 \\
        \leq                       & - \lambda S_{\max} c_u \min_{a,s}\piref(a | s) (\min_{a,s}\pi_{n,1}(a | s))^2 \norm{\frac{d^{\pi^*}_u}{d^{\pi_{n,1}}_u}}^{-1}_\infty \delta_n \\
        \leq                       & -(1-\gamma)\lambda S_{\max} c_u c_{\piref} (c_{\pitheta}^F)^2 \norm{\frac{d^{\pi^*}_u}{u}}^{-1}_\infty \delta_n.
    \end{align*}
    The above inequality can be arranged as
    \begin{align*}
        \delta_{n+1}\leq & \left(1-  (1-\gamma)\lambda S_{\max} c_u c_{\piref} (c_{\pitheta}^F)^2 \norm{\frac{d^{\pi^*}_u}{u}}^{-1}_\infty\right)\delta_n   \\
        \leq             & \exp\left(-(1-\gamma)\lambda S_{\max} c_u c_{\piref} (c_{\pitheta}^F)^2 \norm{\frac{d^{\pi^*}_u}{u}}^{-1}_\infty\right)\delta_n  \\
        \leq             & \exp\left(-(n-1)(1-\gamma)\lambda S_{\max} c_u c_{\piref} (c_{\pitheta}^F)^2 \norm{\frac{d^{\pi^*}_u}{u}}^{-1}_\infty\right)\delta_1.
    \end{align*}
    By converting the initial distribution from \(u\) to \(\rho\), we obtain the final result.
\end{proof}

\subsection{PPO-Clip for reverse KL}
When the regularizer is reverse KL-divergence, the surrogate function follows the form 
\begin{align*}
     \mathcal{L}_{n}(\theta) &:= \E_{\pi_{n,1}}  \Big[\min\big\{r_\theta^n(s,a) \tilde{A}^{\pi_{n,1}}_\lambda(s,a), \\
&\qquad \text{clip}\big(r_\theta^n(s,a);1-\varepsilon_l,\blue{1}+\varepsilon_h\big) \tilde{A}^{\pi_{n,1}}_\lambda(s,a)\big\} \\
&\qquad -  \lambda  r_{\theta}^n(s,a) \log \frac{\pi_{\theta}(a | s)}{\piref(a|s)} \Big].
\end{align*}
Different from forward KL, reverse KL remains finite at the simplex boundary, so objective monotonicity alone cannot guarantee a uniform lower bound on action probabilities. Existing exact-gradient analyses obtain global convergence by exploiting trajectory-wise interiority or local regularity \cite{agarwal2021pg,mei2020global,ding2025twophase}, but these arguments do not directly apply to multi-step PPO-Clip. After the first inner update, \(g_{n,k}\) is the gradient of a clipped surrogate at \(\pi_{n,1}\), rather than the exact policy gradient at \(\pi_{n,k}\). We therefore establish trajectory-wide interiority directly for the deterministic full-surrogate gradient by combining a policy-independent surrogate-gradient error bound with a coordinatewise reverse KL boundary-repulsion property. Together with Lemma~\ref{lem:PL-f}, this yields a uniform \L{}ojasiewicz constant along the trajectory and hence provides global linear convergence.

Since \(D_{\rm KL}(\pi(s),\piref(s))\leq L_{\rm ref}\) for every policy \(\pi\), the reverse KL-regularized advantage satisfies \(\left|\tilde A_\lambda^\pi(s,a)\right|\leq\frac{r_{\rm max}+\lambda L_{\rm ref}}{1-\gamma}\) for all \(\pi,s,a\). We therefore use \(\widetilde A_{\max}^R:=\frac{r_{\rm max}+\lambda L_{\rm ref}}{1-\gamma}\) as a uniform upper bound on the absolute value of the reverse KL-regularized advantage.

\begin{lemma}[Coordinatewise reverse KL boundary repulsion]
\label{lem:rkl-boundary-repulsion}
Suppose Assumption~\ref{asup:min-piref} holds. Then, for every inner iterate \(\theta_{n,k}\) generated by \cref{alg:PPO-Clip}, for every state-action pair \((s,a)\) the gradient
\(g_{n,k}:=\nabla_{\theta_{n,k}}\mathcal{L}_n(\theta_{n,k})\)
satisfies \(|g_{n,k}(s,a)|\leq G_\infty\), where
\(G_\infty:=\big(2\widetilde A_{\max}^R+\lambda(e^{-1}+2L_{\rm ref})\big)/(1-\gamma)\),
and \(g_{n,k}(s,a)\geq0\) whenever
\(\pi_{n,k}(a|s)\leq\tau_R\), where
\(\tau_R:=c_{\piref}\exp\big(-2\widetilde A_{\max}^R/\lambda\big)\).
\end{lemma}

\begin{proof}
Fix an inner-loop point \((n,k)\) and a state \(s\), denote
\(p:=\pi_{n,k}(s)\), \(p^0:=\pi_{n,1}(s)\),
\(q:=\piref(s)\), and \(\bH(p):=\diag(p)-pp^\top\).
For each action, define the clipped coefficient \(b_a\) by
\[
b_a:=\left.\partial_r\min\!\left\{r\tilde A_\lambda^{\pi_{n,1}}(s,a),\operatorname{clip}\big(r;1-\varepsilon_l,1+\varepsilon_h\big)\tilde A_\lambda^{\pi_{n,1}}(s,a)\right\}\right|_{r=p_a/p_a^0} \in[-\widetilde A_{\max}^R,\widetilde A_{\max}^R].
\]
The old-policy probability \(p_a^0\) in the sampling measure cancels with the importance ratio, so the exact state block of the surrogate gradient is
\begin{equation}
g_{n,k}(s,\cdot) =\frac{d_u^{\pi_{n,1}}(s)}{1-\gamma} \bH(p)\left(b-\lambda\log\frac{p}{q}\right), \label{eq:rkl-exact-state-block}
\end{equation}
where the logarithm and division are coordinatewise. Since \([\bH(p)v]_a=p_a(v_a-p^\top v)\), \eqref{eq:rkl-exact-state-block} gives
\begin{equation}
g_{n,k}(s,a) =\frac{d_u^{\pi_{n,1}}(s)}{1-\gamma}p_a \left[ b_a-\sum_{a'\in\cA}p_{a'}b_{a'} -\lambda\left(\log\frac{p_a}{q_a}-\sum_{a'\in\cA}p_{a'}\log\frac{p_{a'}}{q_{a'}}\right) \right]. \label{eq:rkl-coordinate-gradient}
\end{equation}
If \(p_a\leq\tau_R\), then \(q_a\geq c_{\piref}\),
\(\sum_{a'\in\cA}p_{a'}\log(p_{a'}/q_{a'})\geq0\), and
\[
 b_a-\sum_{a'\in\cA}p_{a'}b_{a'} -\lambda\left(\log\frac{p_a}{q_a}-\sum_{a'\in\cA}p_{a'}\log\frac{p_{a'}}{q_{a'}}\right) \geq-2\widetilde A_{\max}^R-\lambda\log\frac{p_a}{q_a} \geq0.
\]
Therefore, \(g_{n,k}(s,a)\geq0\) whenever \(p_a\leq\tau_R\).
It remains to bound a coordinate before it enters the low-probability region. Equation~\eqref{eq:rkl-coordinate-gradient} and the estimates \(p_a|\log p_a|\leq e^{-1}\), \(p_a|\log q_a|\leq L_{\rm ref}\), and \(\sum_{a'\in\cA}p_{a'}\log(p_{a'}/q_{a'})\leq L_{\rm ref}\) imply
\begin{align}
|g_{n,k}(s,a)| &\leq\frac{1}{1-\gamma} \left( 2\widetilde A_{\max}^R p_a +\lambda p_a\left|\log\frac{p_a}{q_a}\right| +\lambda p_a\sum_{a'\in\cA}p_{a'}\log\frac{p_{a'}}{q_{a'}} \right) \leq G_\infty. \label{eq:rkl-coordinate-bound}
\end{align}
Both bounds are policy independent.
\end{proof}

For every inner iterate \((n,k)\), define the centered logits by
\[
\widehat\theta_{n,k}(s,a) :=\theta_{n,k}(s,a) -\frac{1}{\abs{\cA}}\sum_{b\in\cA}\theta_{n,k}(s,b), \qquad \forall s\in\cS,\ a\in\cA.
\]
By the shift invariance of the softmax parameterization, the policies \(\pi_{\widehat\theta_{n,k}}\) and \(\pi_{\theta_{n,k}}\) are identical. Such a property provides a way for us to bound the policy through the centered policy.

\begin{lemma}[Uniform interiority of the policy trajectory]
\label{lem:rkl-uniform-interiority}
Under the setting of Lemma~\ref{lem:rkl-boundary-repulsion}, suppose \(\eta_{n,k}\geq0\) and \(\eta_{n,k}\leq\bar\eta\) for every inner update. Set \(\underline\theta(\bar\eta):=\min\{0,\allowbreak\min_{s,a}\widehat\theta_{1,1}(s,a),\allowbreak\log(\abs{\cA}\tau_R)-\bar\eta G_\infty\}\).
Then the centered logits satisfy
\(\underline\theta(\bar\eta) \leq\widehat\theta_{n,k}(s,a) \leq-(\abs{\cA}-1)\underline\theta(\bar\eta)\),
and the complete trajectory is uniformly interior with
\begin{equation}
\inf_{n,k,s,a}\pi_{n,k}(a|s) \geq c_{\min}(\bar\eta) :=\frac{1}{1+(\abs{\cA}-1)\exp\big(-\abs{\cA}\underline\theta(\bar\eta)\big)}>0, \label{eq:rkl-uniform-interiority}
\end{equation}
\end{lemma}

\begin{proof}
Equation~\eqref{eq:rkl-exact-state-block} implies
\(\sum_{a\in\cA}g_{n,k}(s,a)=0\). Hence centering is preserved by every update as
\[\widehat\theta_{n,k+1}(s,a) =\widehat\theta_{n,k}(s,a) +\eta_{n,k}g_{n,k}(s,a), \quad \sum_{a\in\cA}\widehat\theta_{n,k}(s,a)=0.\]
We prove the claimed centered-logit bounds by induction over all inner updates. The lower bound holds at initialization by definition. \emph{If \(p_a\leq\tau_R\)}, the boundary-repulsion property in Lemma~\ref{lem:rkl-boundary-repulsion} shows that \(g_{n,k}(s,a)\geq0\), so the corresponding centered logit does not decrease. \emph{If \(p_a>\tau_R\)}, Jensen's inequality and the zero-sum condition give \(\sum_{b\in\cA}\exp(\widehat\theta_{n,k}(s,b))\geq \abs{\cA}\), and hence 
\[\widehat\theta_{n,k}(s,a)=\log p_a+\log\sum_{b\in\cA}\exp(\widehat\theta_{n,k}(s,b))\geq\log(\abs{\cA}\tau_R).\]
Since \(\eta_{n,k}\leq\bar\eta\), \eqref{eq:rkl-coordinate-bound} gives
\[
\widehat\theta_{n,k+1}(s,a)\geq\log(\abs{\cA}\tau_R)-\bar\eta G_\infty\geq\underline\theta(\bar\eta).
\]
This proves the lower bound at the next iterate. The upper bound follows from the zero-sum identity by summing the lower bounds for the other \(\abs{\cA}-1\) coordinates. Finally,
\[
\pi_{n,k}(a|s) =\frac{1}{1+\sum_{b\neq a} \exp(\widehat\theta_{n,k}(s,b)-\widehat\theta_{n,k}(s,a))} \geq \frac{1}{1+(\abs{\cA}-1)e^{-\abs{\cA}\underline\theta(\bar\eta)}} =c_{\min}(\bar\eta),
\]
which proves the result.
\end{proof}

\begin{theorem}[Global convergence under reverse KL]
\label{thm:global-rkl}
Suppose Assumptions~\ref{asup:min-piref} and~\ref{asup:min-u} hold and \(0<\varepsilon_l<1\). Consider \cref{alg:PPO-Clip} initialized at any finite \(\theta_{1,1}\in\R^{\dS\dA}\). There exists an explicit constant \(\bar S_R>0\), depending only on the problem and algorithm parameters, such that if \(S:=\sum_{k=1}^K\eta_{n,k}\in(0,\bar S_R]\), then the complete policy trajectory is uniformly interior as 
\(\inf_{n,k,s,a}\pi_{n,k}(a|s)\geq c_{\min}(S)>0.\)
Moreover, defining
\(\beta := (1-\gamma)\lambda c_u c_{\min}(S)^2 \left\|\frac{d_u^{\pi^*}}{u}\right\|_\infty^{-1},\)
for every initial distribution \(\rho\), the optimality gap decreases as 
\begin{align}
\tilde V_\lambda^{\pi^*}(\rho) -\tilde V_\lambda^{\pi_{n,1}}(\rho) \leq \frac{ \tilde V_\lambda^{\pi^*}(u) -\tilde V_\lambda^{\pi_{1,1}}(u) }{ c_u(1-\gamma) } e^{-\beta S(n-1)}. \label{eq:rkl-global-rate}
\end{align}
And the policy distance satisfies
\begin{align}
\sum_s \left\| \pi_{n,1}(s)-\pi^*(s) \right\|_2^2 \leq \frac{ 2\big( \tilde V_\lambda^{\pi^*}(u) -\tilde V_\lambda^{\pi_{1,1}}(u) \big) }{ \lambda c_u } e^{-\beta S(n-1)}. \label{eq:rkl-policy-rate}
\end{align}
Hence, \(\pi_{n,1}\) converges linearly to the unique optimal policy \(\pi^*\).
\end{theorem}

\begin{proof}
    For \(n\geq1\), define
\(\delta_n^u := \tilde V_\lambda^{\pi^*}(u) - \tilde V_\lambda^{\pi_{n,1}}(u).\)
Let
\(C_g := \widetilde A_{\max}^R (\sqrt{2}\kappa_\varepsilon+3) + \lambda \left( 6L_{\rm ref} +3\log|\cA| +2 \right),\)
and set
\(\bar S_R := \min\left\{ \frac{1-\gamma}{8C_g}, \frac{1}{4L_R} \right\}.\)
We first control the discrepancy between the PPO-Clip surrogate
gradient and the exact policy gradient at the outer iterate by
separating the clipping, policy-shift, and regularization errors.
Let \(\Delta_{n,k}:=\theta_{n,k}-\theta_{n,1}\). In the decomposition
\[
g_{n,k}-\nabla_\theta\tilde V_\lambda^{\pi_{n,1}}(u) =X_{n,k}+Y_{n,k}+\widetilde Z_{n,k},
\]
the terms \(X_{n,k}\) and \(Y_{n,k}\) represent the clipping and policy-shift errors, respectively. They have the same form as the corresponding terms in the proof of Lemma~\ref{lem:descent-lemma}, so their bounds follow from the same arguments, with \(\widetilde A_{\max}^F\) replaced by \(\widetilde A_{\max}^R\):
\[
\norm{X_{n,k}}_2 \leq\frac{\sqrt{2}\widetilde A_{\max}^R}{1-\gamma} \kappa_\varepsilon \norm{\Delta_{n,k}}_2, \qquad \norm{Y_{n,k}}_2 \leq\frac{3\widetilde A_{\max}^R}{1-\gamma} \norm{\Delta_{n,k}}_2.
\]
The remaining term \(\widetilde Z_{n,k}\) is the regularization error, whose form depends on the regularizer. For reverse KL, the regularization term is separable across states. For a fixed state \(s\), write the reverse KL regularizer as the following function of the corresponding statewise logits:
\[
\mathcal R_s(\theta_s) :=D_{\rm KL}\big(\pi_\theta(s),\piref(s)\big) =\sum_a\pi_\theta(a|s) \log\frac{\pi_\theta(a|s)}{\piref(a|s)}.
\]
Since \(\sum_a\nabla_{\theta_s}\pi_\theta(a|s)=0\),
\[
\nabla_{\theta_s}\mathcal R_s(\theta_s) =\sum_a\nabla_{\theta_s}\pi_\theta(a|s) \log\frac{\pi_\theta(a|s)}{\piref(a|s)}.
\]
Consequently, the \(s\)-block of \(\widetilde Z_{n,k}\) is
\begin{align*}
\widetilde Z_{n,k}(s,\cdot) =-\frac{\lambda}{1-\gamma}d_u^{\pi_{n,1}}(s) \big[ \nabla\mathcal R_s(\theta_{n,k}(s,\cdot)) -\nabla\mathcal R_s(\theta_{n,1}(s,\cdot)) \big].
\end{align*}
By \eqref{smo:rkl-3}, the Hessian of this statewise regularization function satisfies \(\norm{\nabla^2\mathcal R_s(\theta_s)}_2\leq 6L_{\rm ref}+3\log\abs{\cA}+2\) globally. Applying the mean-value theorem to the difference of gradients in each state block, using \(d_u^{\pi_{n,1}}(s)\leq1\), and then summing the squared block norms gives
\begin{equation}
\norm{\widetilde Z_{n,k}}_2 \leq \frac{\lambda(6L_{\rm ref}+3\log\abs{\cA}+2)}{1-\gamma} \norm{\Delta_{n,k}}_2. \label{eq:rkl-global-z}
\end{equation}
Combining the preceding bounds with \eqref{eq:rkl-global-z} gives
\begin{equation}
\norm{g_{n,k}-\nabla_\theta\tilde V_\lambda^{\pi_{n,1}}(u)}_2 \leq \frac{C_g}{1-\gamma} \norm{\Delta_{n,k}}_2. \label{eq:rkl-global-inexactness}
\end{equation}
Since \(S\leq(1-\gamma)/(8C_g)\), the same inner-loop induction as in Lemma~\ref{lem:descent-lemma} yields \(\norm{\Delta_{n,k}}_2\leq2S\norm{\nabla_\theta\tilde V_\lambda^{\pi_{n,1}}(u)}_2\).
Substitution into \eqref{eq:rkl-global-inexactness} gives the uniform relative-error estimate
\begin{equation}
\norm{g_{n,k}-\nabla_\theta\tilde V_\lambda^{\pi_{n,1}}(u)}_2 \leq\frac14\norm{\nabla_\theta\tilde V_\lambda^{\pi_{n,1}}(u)}_2. \label{eq:rkl-relative-error}
\end{equation}
Using the \(L_R\)-smoothness of \(\tilde V_\lambda\) and the update \(\theta_{n+1,1}-\theta_{n,1}=\sum_k\eta_{n,k}g_{n,k}\), we have
\begin{align}
\delta_{n+1}^u-\delta_n^u &\leq -\left\langle \nabla_\theta\tilde V_\lambda^{\pi_{n,1}}(u), \sum_k\eta_{n,k}g_{n,k} \right\rangle +\frac{L_R}{2} \left\|\sum_k\eta_{n,k}g_{n,k}\right\|_2^2 \notag\\
&\leq -\frac{3S}{4}\|\nabla_\theta\tilde V_\lambda^{\pi_{n,1}}(u)\|_2^2 +\frac{25L_RS^2}{32} \|\nabla_\theta\tilde V_\lambda^{\pi_{n,1}}(u)\|_2^2 \notag\\
&\leq -\frac{S}{2} \|\nabla_\theta\tilde V_\lambda^{\pi_{n,1}}(u)\|_2^2. \label{eq:rkl-gap-descent}
\end{align}
Here the second inequality uses \eqref{eq:rkl-relative-error}, and the third uses \(S\leq1/(4L_R)\).
Thus, the objective gap decreases at every outer iteration under a fixed, policy-independent stepsize.

We next use the coordinatewise sign property in
Lemma~\ref{lem:rkl-boundary-repulsion} to obtain a uniform lower
bound on all action probabilities along the entire trajectory.
Because the inner stepsizes are nonnegative and sum to \(S\),
each satisfies \(\eta_{n,k}\leq S\). Applying
Lemmas~\ref{lem:rkl-boundary-repulsion} and
\ref{lem:rkl-uniform-interiority} with \(\bar\eta=S\) gives
\[
\inf_{n,k,s,a}\pi_{n,k}(a|s)\geq c_{\min}(S)>0.
\]
With this trajectory-wide interiority bound, the non-uniform
\L{}ojasiewicz inequality in Lemma~\ref{lem:PL-f} becomes uniform
along the outer iterates and yields a geometric contraction of the
value gap.
For reverse KL, the strong-concavity coefficient in Lemma~\ref{lem:PL-f} is \(c_m=1\). Applying that lemma with \(\rho=u\), together with Lemma~\ref{lem:rkl-uniform-interiority} at \(\bar\eta=S\), so that \(\min_{s,a}\pi_{n,1}(a|s)\geq c_{\min}(S)\), yields
\[
\norm{\nabla_\theta\tilde V_\lambda^{\pi_{n,1}}(u)}_2^2 \geq 2\lambda c_u c_{\min}(S)^2 \norm{\frac{d_u^{\pi^*}}{d_u^{\pi_{n,1}}}}_\infty^{-1} \delta_n^u.
\]
Moreover, \(d_u^{\pi_{n,1}}(s)\geq(1-\gamma)u(s)\), so
\[
\norm{\frac{d_u^{\pi^*}}{d_u^{\pi_{n,1}}}}_\infty^{-1} \geq(1-\gamma)\norm{\frac{d_u^{\pi^*}}{u}}_\infty^{-1}.
\]
By the definition of \(\beta\) in the theorem, we obtain the trajectory-wide PL inequality
\(\norm{\nabla_\theta\tilde V_\lambda^{\pi_{n,1}}(u)}_2^2\geq2\beta\delta_n^u\).
Combining the above PL inequality with \eqref{eq:rkl-gap-descent} gives
\[\delta_{n+1}^u\leq(1-\beta S)\delta_n^u \leq e^{-\beta S}\delta_n^u.\]
If \(\delta_n^u>0\), the nonnegativity of \(\delta_{n+1}^u\) in the first inequality implies \(\beta S\leq1\); if \(\delta_n^u=0\), the desired conclusion is already exact. Thus, in either case,
\(\delta_n^u\leq\delta_1^u\exp\big(-\beta S(n-1)\big)\).

Finally, we convert the value-gap contraction into a
policy-distance bound through the statewise regularized Bellman
objective, and then transfer the value-gap estimate from \(u\)
to an arbitrary initial distribution \(\rho\).
Define the statewise regularized Bellman objective
\[
\mathcal B_s(p) :=\sum_a p_a\left(r(s,a)+\gamma\sum_{s'}\cP(s'\mid s,a) \tilde V_\lambda^{\pi^*}(s')\right) -\lambda D_{\rm KL}\big(p,\piref(s)\big).
\]
By the regularized Bellman optimality equation,
$\tilde V_\lambda^{\pi^*}(s)
=\max_{p\in\Delta(\cA)}\mathcal B_s(p)$,
and hence \(\pi^*(s)\) is a maximizer of \(\mathcal B_s\).
Moreover, \(\mathcal B_s\) is \(\lambda\)-strongly concave,
so this maximizer is unique and
\[
\mathcal B_s(\pi^*(s))-\mathcal B_s(p) \geq \frac{\lambda}{2}\norm{p-\pi^*(s)}_2^2.
\]
Since \(\mathcal B_s(\pi^*(s)) =\tilde V_\lambda^{\pi^*}(s)\), we obtain
\[
\tilde V_\lambda^{\pi^*}(s)-\mathcal B_s(p) \geq \frac{\lambda}{2}\norm{p-\pi^*(s)}_2^2.
\]
Using the performance-difference identity with initial distribution \(u\), we obtain
\begin{align*}
\delta_n^u =\frac{1}{1-\gamma}\sum_s d_u^{\pi_{n,1}}(s) \left[\tilde V_\lambda^{\pi^*}(s) -\mathcal B_s\big(\pi_{n,1}(s)\big)\right] \geq\frac{\lambda c_u}{2} \sum_s\norm{\pi_{n,1}(s)-\pi^*(s)}_2^2.
\end{align*}
Here the inequality uses the strong concavity of \(\mathcal B_s\) and
\(d_u^{\pi_{n,1}}(s)\geq(1-\gamma)u(s)\geq(1-\gamma)c_u\).
Combining the last two displays proves \eqref{eq:rkl-policy-rate}. For every policy \(\pi\), the same performance-difference argument as in \eqref{eq:change-initial-distribution} gives
\[
\tilde V_\lambda^{\pi^*}(\rho)-\tilde V_\lambda^\pi(\rho) \leq\frac{1}{c_u(1-\gamma)} \left[ \tilde V_\lambda^{\pi^*}(u)-\tilde V_\lambda^\pi(u) \right].
\]
Iterating the recursion and applying this distribution-conversion inequality proves \eqref{eq:rkl-global-rate}.
\end{proof}

\begin{remark}[Initialization at the reference policy]
\label{rmk:rkl-reference-initialization}
Under the initialization \(\pi_{1,1}=\piref\) in \cref{alg:PPO-Clip}, the centered logits can be chosen as
\(\widehat\theta_{1,1}(s,a) = \log\piref(a|s) -\frac{1}{\abs{\cA}}\sum_b\log\piref(b|s).\)
Since the average of \(\log\piref(b|s)\) is nonpositive, it can be verified 
\(\min_{s,a}\widehat\theta_{1,1}(s,a)\geq \log c_{\piref}.\)
So \(c_{\min}(S)\) in \cref{thm:global-rkl} can be explicitly bounded in terms of the problem parameters and the fixed stepsize \(S\). 
\end{remark}

The fixed-stepsize analysis above uses a policy-independent worst-case bound on the surrogate-gradient error, which leads to a uniform but potentially conservative stepsize restriction. To obtain a less conservative rule, we next allow the total inner-loop stepsize to vary across outer iterations according to the current policy. Specifically, the adaptive quantity is
\(S_n:=\sum_{k=1}^K\eta_{n,k},\)
which will be chosen based on an iterate-dependent bound on the surrogate-gradient error. For this purpose, we retain the dependence of this error bound on the current outer policy and define
\[
C_{g,n} :=\widetilde A_{\max}^R(\sqrt{2}\kappa_\varepsilon+3) +3\lambda\max_{s,a} \left| \log\frac{\piref(a|s)}{\pi_{n,1}(a|s)} \right| +2\lambda.
\]
The quantity \(C_{g,n}\) measures the local surrogate-gradient error at outer iteration \(n\) and will be used to select \(S_n\). We will show that \(C_{g,n}\) remains uniformly bounded along the trajectory, so the adaptive stepsizes are also uniformly bounded away from zero.

\begin{corollary}[Global  convergence with adaptive stepsizes]
\label{cor:rkl-adaptive-global}
Under the assumptions of \cref{thm:global-rkl}, consider the
same PPO-Clip scheme with stepsizes \(\eta_{n,k}\) satisfy 
\(S_n:=\sum_{k=1}^K\eta_{n,k} =\min\left\{\frac{1-\gamma}{8C_{g,n}},\frac{1}{4L_R}\right\}\).
Then the generated trajectory satisfies
\(\pi_{n,k}(a|s)\geq\bar c:=c_{\min}\!\left(\frac{1}{4L_R}\right)\). Moreover, there exists a constant
\(S_{\min}>0\),  such that
\(S_n\geq S_{\min}\).
Consequently, both the value gap and the squared policy distance
converge globally at a linear rate.
\end{corollary}

\begin{proof}
With the adaptive rule, the reverse KL regularization error satisfies
\[
\norm{\widetilde Z_{n,k}}_2 \leq \frac{\lambda}{1-\gamma} \left( 3\max_{s,a} \left| \log\frac{\piref(a|s)}{\pi_{n,1}(a|s)} \right| +2 \right) \norm{\Delta_{n,k}}_2.
\]
Together with the bounds for \(X_{n,k}\) and \(Y_{n,k}\) used in
the proof of \cref{thm:global-rkl}, the inexact-gradient estimate
\eqref{eq:rkl-global-inexactness} remains valid with \(C_g\) replaced
by \(C_{g,n}\). Since \(S_n\leq(1-\gamma)/(8C_{g,n})\), the inner-loop
induction from Lemma~\ref{lem:descent-lemma} then yields
\[
\norm{ g_{n,k} -\nabla_\theta\tilde V_\lambda^{\pi_{n,1}}(u) }_2 \leq \frac14 \norm{ \nabla_\theta \tilde V_\lambda^{\pi_{n,1}}(u) }_2.
\]
Together with \(S_n\leq1/(4L_R)\), the smoothness argument
leading to \eqref{eq:rkl-gap-descent} gives
\[
\delta_{n+1}^u-\delta_n^u \leq -\frac{S_n}{2} \norm{ \nabla_\theta \tilde V_\lambda^{\pi_{n,1}}(u) }_2^2.
\]
As \(S_n\leq1/(4L_R)\), every inner stepsize satisfies
\(\eta_{n,k}\leq1/(4L_R)\).
Lemma~\ref{lem:rkl-uniform-interiority} therefore gives
\(\pi_{n,k}(a|s) \geq c_{\min}\!\left(\frac{1}{4L_R}\right) =:\bar c>0.\)
Consequently,
\[
\max_{s,a} \left| \log\frac{\piref(a|s)}{\pi_{n,1}(a|s)} \right| \leq L_{\rm ref}+\log\frac1{\bar c},
\]
and hence
\[
C_{g,n} \leq \widetilde A_{\max}^R(\sqrt{2}\kappa_\varepsilon+3) +3\lambda \left(L_{\rm ref}+\log\frac{1}{\bar c}\right) +2\lambda =:C_{g,\max}.
\]
By Lemma~\ref{lem:PL-f} and the uniform bound
\(\pi_{n,k}(a|s)\geq\bar c\), the PL inequality reads
\(\norm{\nabla_\theta\tilde V_\lambda^{\pi_{n,1}}(u)}_2^2\allowbreak \geq2\bar\beta\,\delta_n^u\), where
\(\bar\beta:=(1-\gamma)\lambda c_u\bar c^2 \norm{\frac{d_u^{\pi^*}}{u}}_\infty^{-1}\).
Combining the last two inequalities yields
\(\delta_{n+1}^u\leq(1-\bar\beta S_n)\delta_n^u \leq e^{-\bar\beta S_n}\delta_n^u\).
Therefore,
\[
\delta_n^u \leq \delta_1^u \exp\Big( -\bar\beta\sum_{j=1}^{n-1}S_j \Big) \leq \delta_1^u e^{-\bar\beta S_{\min}(n-1)}.
\]
The distribution-conversion argument in the proof of
\cref{thm:global-rkl} then gives
\[
\tilde V_\lambda^{\pi^*}(\rho) -\tilde V_\lambda^{\pi_{n,1}}(\rho) \leq \frac{\delta_1^u}{c_u(1-\gamma)} e^{-\bar\beta S_{\min}(n-1)},
\]
while the statewise strong-concavity argument yields
\[
\sum_s \norm{\pi_{n,1}(s)-\pi^*(s)}_2^2 \leq \frac{2\delta_1^u}{\lambda c_u} e^{-\bar\beta S_{\min}(n-1)}. \qedhere
\]
\end{proof}

\section{Conclusion}
We study the theoretical properties of PPO-Clip in infinite-horizon discounted MDPs with \(f\)-divergence regularization under the softmax policy parameterization. Motivated by the use of clipping and divergence regularization in modern policy optimization we develop our analysis in a general RL setting.  We formulate the \(f\)-divergence regularized value function and establish its structural properties, including a non-uniform Lipschitz smoothness condition and a \L{}ojasiewicz inequality under mild conditions on the \(f\)-function.

Building on these properties, we analyze the global convergence of a deterministic actor-only PPO-Clip algorithm under forward and reverse KL regularization. For forward KL, we establish an \(O(1/N)\) stationarity rate and a global linear convergence of the value gap. For reverse KL, we identify a coordinatewise boundary-repulsion property that keeps the policy trajectory uniformly away from the boundary, which leads to global linear convergence in both the value gap and the squared policy distance from any finite softmax initialization. These results provide a theoretical characterization of how clipping and divergence regularization interact with the geometry of softmax policy optimization.

\end{document}